\documentclass[11pt,dvips]{amsart}
\usepackage{amsmath}
\usepackage{amssymb}
\usepackage{amscd}
\usepackage{palatino}
\usepackage{euler}
\usepackage{latexsym}
\usepackage{epsfig}
\usepackage{graphics}
\usepackage{amsfonts}
\usepackage{psfrag}

\input xy
\xyoption{all}

\numberwithin{equation}{section}
\numberwithin{figure}{section}

\newtheorem{theorem}{Theorem}[section]
\newtheorem{lemma}[theorem]{Lemma}
\newtheorem{proposition}[theorem]{Proposition}
\newtheorem{corollary}[theorem]{Corollary}

\newtheorem{definition}[theorem]{Definition}

\newtheorem{example}[theorem]{Example}
\newtheorem{remark}[theorem]{Remark}
\newcommand\dfn{\bf}

\newcommand{\C}{{\mathbb{C}}}

\newcommand{\Z}{{\mathbb{Z}}}
\newcommand{\Q}{{\mathbb{Q}}}
\newcommand{\R}{{\mathbb{R}}}

\newcommand{\T}{{\mathbb{T}}}
\renewcommand{\P}{{\mathbb{P}}}

\newcommand{\algt}{\mathfrak{t}}

\newcommand{\algk}{\mathfrak{k}}

\newcommand\tensor{\otimes}

\newcommand{\grad}{\rm grad\mbox{ }}

\newcommand{\iso}{\cong}
\newcommand{\into}{\hookrightarrow}
\newcommand{\onto}{\twoheadrightarrow}
\newcommand\lie[1]{{\mathfrak #1}}

\newcommand{\stab}{\mathrm{Stab}}

\newcommand{\age}{\mathrm{age}}
\newcommand{\finitepreorb}{NH_T^{*,\Gamma}}
\newcommand{\rous}{NH_{T}^{*,\diamond}}
\newcommand{\krous}{NH_{K}^{*,\diamond}}
\newcommand{\srous}{NH_{S^1}^{*,\diamond}}
\newcommand{\kfinitepreorb}{NH_K^{*,\Gamma}}
\newcommand{\Trous}{NH_{\T}^{*,\diamond}}

\newcommand\rousT{H^*_T(Y^T) \tensor_{\Z} \Z[T]}
\newcommand{\Year}{{a}}

\newcommand{\finitepiece}{\mbox{$\Gamma$-subring}}
\newcommand{\level}{\ensuremath{\Phi^{-1}(0)}}
\newcommand{\gpiece}{NH_T^{*,g}}
\newcommand{\locfree}{\mathcal{Z}}

\newcommand{\I}{\mathcal{I}}
\newcommand{\J}{\mathcal{J}}

\newcommand\AGV{{$\aleph$}\!\! GV}
\newcommand\To{\longrightarrow}
\newcommand{\RoUS}{\rous}
\newcommand\RoUST{\rousT}
\newcommand\hatT{\widehat T}
\newcommand\Hom{{\rm Hom}}
\newcommand\Span{{\rm Span}}

\newcommand{\todo}[2][0.8]{\vspace{1 mm}\par \noindent
\marginpar{\textsc{To Do}} \framebox{\begin{minipage}[c]{#1
\textwidth} \tt #2 \end{minipage}}\vspace{1 mm}\par}

\begin{document}

\title{Orbifold cohomology of torus quotients}
\author{Rebecca Goldin}
\address{Mathematical Sciences MS 3F2, George Mason University,
  4400 Univ.\ Dr., Fairfax, VA 22030, USA}
\email{rgoldin@math.gmu.edu}
\thanks{RG was  supported by
  NSF-DMS Grant 0305128 and by MSRI}
\author{Tara S.\ Holm}
\address{Department of Mathematics, Cornell University,
  Ithaca, NY 14853-4201, USA}
\email{tsh@math.cornell.edu}
\thanks{TSH was supported in part by a National
  Science Foundation Postdoctoral Fellowship.}
\author{Allen Knutson}
\address{Department of Mathematics, University of California,
  San Diego, CA 92093, USA}
\thanks{AK was supported by NSF grant 0303523.}
\email{allenk@math.ucsd.edu}
\date{\today}

%{\it MSC 2000 Subject Classification}:
%Primary: xxxxx \hspace{0.1in} Secondary: xxxxx, xxxxx
%\newline \mbox{~~~~} {\it Keywords:} equivariant cohomology,
%symplectic reduction, orbifold}\\

%%%%%%%%%%%%%%%%%%%%
% Disclamer
%%%%%%%%%%%%%%%%%%%%
%\begin{center}
%\framebox{
%{\Large\bf DRAFT (\today): DO NOT DISTRIBUTE.}}
%\end{center}

%%%%%%%%%%%%%%%%%%%%%
%  Abstract
%%%%%%%%%%%%%%%%%%%%%

\begin{abstract}
  We introduce the {\dfn inertial cohomology ring} $\rous(Y)$ of a
  stably almost complex manifold carrying an action of a torus $T$.
%  The bigrading on this ring is by real numbers in the first entry,
%  and by elements of $T$ in the second.
  We show that in the case that $Y$ has a locally free action by $T$,
  the inertial cohomology ring is isomorphic to the Chen-Ruan orbifold
  cohomology ring $H_{CR}^*(Y/T)$ (as defined in [Chen-Ruan]) of the
  quotient orbifold $Y/T$.

  For $Y$ a compact Hamiltonian $T$-space, we extend to orbifold
  cohomology two techniques that are standard in ordinary cohomology.
  We show that $\rous(Y)$ has
  a natural ring surjection onto $H_{CR}^*(Y/\!/T)$, where $Y/\!/T$
  is the symplectic reduction of $Y$ by $T$ at a regular value of the
  moment map.  We extend to $\rous(Y)$ the graphical GKM calculus (as
  detailed in e.g. [Harada-Henriques-Holm]), and the kernel
  computations of [Tolman-Weitsman, Goldin].

  We detail this technology in two examples: toric orbifolds and
  weight varieties, which are symplectic reductions of flag manifolds.
  The Chen-Ruan ring has been computed for toric orbifolds, with
  $\Q$ coefficients, in
  [Borisov-Chen-Smith]); we reproduce their results over $\Q$ for all
  symplectic toric orbifolds obtained by reduction by a connected
  torus (though with different computational methods), and extend them
  to $\Z$ coefficients in certain cases, including weighted projective spaces.
\end{abstract}

\maketitle

\setcounter{tocdepth}{1}
\tableofcontents

%% 1. Introduction -- Rebecca will write a sketchy draft
%%    intro.tex
%%        * Including the nicest most computable definition

\section{Introduction}
\label{sec:intro}
In \cite{CR:orbH}, Chen and Ruan introduced orbifold
cohomology groups along with a product structure as part of a program
to understand orbifold string theory. This {\dfn Chen-Ruan orbifold cohomology ring\footnote{Here,
we use $H^*_{CR}$ to denote the Chen-Ruan orbifold cohomology ring.
In \cite{CR:orbH}, Chen and Ruan call this ring {\dfn orbifold
  cohomology} and denote it $H^*_{orb}$, but this name and notation
have been used multiply in the literature.}}
$H^*_{CR}(X)$ is the degree $0$ part of the quantum cohomology of
the orbifold $X$ (which when $X$ is a manifold, reduces to the
ordinary cohomology), and as such, one of its subtlest properties is
the associativity of its product.  
%% It is conjectured that
%% $H^*_{CR}(X;\C)$ with complex  
%% coefficients is isomorphic as a ring to the ordinary cohomology
%% $H^*(\widetilde X;\C)$ of a hyperkahler resolution $\widetilde{X}$ of $X$,
%% when one exists.  More generally, for a crepant resolution
%% $\widetilde{X}$ of $X$, Ruan's quantum minimal model  
%% conjecture states that there is a ``quantum-corrected cohomology ring"
%% of $\widetilde{X}$ that is isomorphic to $H^*_{CR}(X)$ (see
%% \cite{R:crepantresolutions,R:stringy}). In this way,
%% $H^*_{CR}(X)$ should record data about
%% some of the simplest kind of singularities of blowdowns, namely
%% orbifold singularities.  
It was originally conjectured that
$H^*_{CR}(X;\C)$ with complex
coefficients is isomorphic as a ring to the ordinary cohomology
$H^*(\widetilde X;\C)$ of a crepant resolution $\widetilde{X}$ of $X$,
when one exists (see e.g.\ \cite{CR:orbH,BCS:toricvarieties}).
In this way,
$H^*_{CR}(X)$ should record data about
some of the simplest kind of singularities of blowdowns, namely
orbifold singularities.
For example, ``simple'' singularities in
codimension $2$ are all orbifold singularities. These arise by
blowing down ADE diagrams of rational curves in a surface.
(Ruan's quantum minimal model conjecture \cite{R:crepantresolutions,R:stringy}
modifies this conjecture slightly, and involves corrections from the
quantum cohomology of $\widetilde{X}$.)

Fantechi and G\"ottsche simplified the presentation of $H^*_{CR}(X)$ in
\cite{FG:globalquotients}, in the case that $X$ is the global
quotient of a complex manifold by a finite (possibly nonabelian)
group. In the algebraic category, Abramovich, Graber and Vistoli
\cite{AGV:orbifoldquantumproducts} described an analogous story to
Chen-Ruan's for Deligne-Mumford stacks. % and moduli spaces of curves.
Borisov, Chen and Smith \cite{BCS:toricvarieties}
used the \AGV\ prescription to describe explicitly the Chen-Ruan
cohomology for toric Deligne-Mumford stacks.

The goal of this paper is to simplify the presentation of the Chen-Ruan
cohomology ring for those orbifolds that occur as a
global quotient by an abelian compact Lie group.
These orbifolds were already intensively studied by Atiyah in
\cite{Atiyah},
where he essentially computed an index theorem for them,
using a Chern character map
taking values in what we now recognize to be their Chen-Ruan cohomology groups
(which did not have a general definition at the time).
%Atiyah does not appear to make reference to the Chen-Ruan product.

Our interest in this family of orbifolds is due to their origin in the
study of symplectic reductions of Hamiltonian $T$-spaces.  Recall that
a symplectic manifold $(Y,\omega)$ carrying an action of a torus $T$
is a {\dfn Hamiltonian $T$-space} if there is an invariant map
$\Phi: Y\rightarrow \mathfrak{t}^*$ from $Y$ to the dual of the Lie algebra
of $T$ satisfying
\begin{equation}\label{eq:momentmapcondition}
  d\langle\Phi ,\xi\rangle = \iota_{V_\xi}\omega
\end{equation}
for all $\xi\in \mathfrak{t}$, where $V_\xi$ is the vector field on
$Y$ generated by $\xi$. Throughout this paper, we will assume that
some component $\langle\Phi ,\xi\rangle$ of the moment map $\Phi$ is
proper and bounded below. We call a Hamiltonian $T$-space whose moment
map satisfies this condition a {\dfn proper Hamiltonian T-space}.
The most important examples are smooth projective varieties $Y$
carrying a linear $T$-action; the symplectic form is the Fubini-Study
form from the ambient projective space,
and properness follows from compactness.

It follows from (\ref{eq:momentmapcondition}) that, for any regular
value $\mu$ of the moment map, $\Phi^{-1}(\mu)$ is a submanifold of
$Y$ with a locally free $T$-action. In particular, any point in the
level set has at most a finite stabilizer in $T$. The symplectic
reduction
$$
Y/\!/T(\mu):= \Phi^{-1}(\mu)/T
$$
is thus an orbifold.
In particular, many (but not all) toric orbifolds may be
obtained by symplectic reduction of manifolds.

Recall that for any $T$-space $Y$, the inclusion
$i:Y^T\hookrightarrow Y$ of the fixed point set induces a map
$$
i^*:H_T^*(Y)\hookrightarrow H_T^*(Y^T)
$$
in equivariant cohomology. If $Y$ has the property that $i^*$ is
{\em injective} (over $\Q$, $\Z$, etc.), then we call $Y$ {\dfn
equivariantly injective} (over $\Q$, $\Z$, etc.).
We say that $Y$ is {\dfn equivariantly formal} (with respect to its
$T$-action) if the $E_2$ term of the Leray-Serre spectral sequence
$$
Y\hookrightarrow Y\times_TET\rightarrow BT
$$
collapses, implying $H_T^*(Y) \iso H^*(Y)\otimes H_T^*(pt)$ as
modules over $H_T^*(pt)$. Over $\Q$, equivariant formality implies
equivariant injectivity (see \cite{GS:supersymmetry}). In particular,
proper Hamiltonian $T$-spaces are always equivariantly injective over
$\Q$.

For any $g\in T$,
let $Y^g$ denote the fixed point set of the $g$ action on $Y$. We
will say that $Y$ is {\dfn robustly equivariantly injective} if, for
every $g\in T$, the $T$-invariant submanifold $Y^g$ is equivariantly
injective. When $Y$ is a proper Hamiltonian space and $H^*(Y^T)$ is
free, $Y$ is equivariantly injective over $\Z$, and even robustly
equivariantly injective (since $Y^g$ is itself a proper Hamiltonian space
for every $g$).

Note that not all equivariantly injective spaces are robustly so: a
counterexample is $S^1$ acting on $\R\P^2$ by ``rotation'' (induced
from the rotation action on the double cover $S^2$). The points
fixed by the element of order $2$ form the set $\R\P^1 \cup \{pt\}$,
which is not equivariantly injective. We thank C. Allday and V.
Puppe for each discovering this example and sharing it with us.

Our main contribution to the study of Chen-Ruan cohomology is the
definition of the {\dfn inertial cohomology\footnote{In the
    announcement \cite{GHK:ober} of these results, we used the term
    {\dfn preorbifold cohomology}.  We believe {\dfn inertial} is more
    suggestive, referring to the inertia orbifold, whose cohomology we
    are studying.  We use the notation $NH^*$ for i{\bf n}ertial
    cohomology because $IH^*$ is the standard notation for
    intersection cohomology.}}  
of a stably almost complex
manifold $Y$, and in particular of a proper Hamiltonian $T$-space, denoted
$\rous(Y)$.
The inertial cohomology is defined as an $H_T^*(pt)$-module
(but not a ring) by
\begin{align*}
\rous(Y) :&= \bigoplus_{g\in T}  H_T^*(Y^g).
\end{align*}
The product structure and even the grading on $\rous(Y)$
are rather complicated, and we
leave their definition to Section~\ref{sec:smileproduct}.
%We shall see that $\rous(Y)$ is an associative, bigraded ring.
The grading $\ast$ is by {\em real} numbers,
and $\diamond$ by elements of $T$.
The summand $H_T^*(Y^g)$ in the above is the
$NH_T^{*,g}(Y) = \bigoplus_{r\in\R} NH_T^{r,g}(Y)$ part of $\rous(Y)$.

The collection of restriction maps $H^*_T(Y^g)\rightarrow H_T^*(Y^T)$
produces a map
\begin{equation}\label{eq:introinjective}
 \rous(Y) \longrightarrow \bigoplus_{g\in T}
H^*_T(Y^T) \iso \rousT.
\end{equation}
In Section~\ref{sec:starproduct} we define a product $\star$ on the
target space.  In the case of $Y$ robustly equivariantly injective,
the map (\ref{eq:introinjective}) is {\em injective}, and $\star$
pulls back to a product on $\rous(Y)$. Most importantly, this product
is easy to compute.

In Section~\ref{sec:smileproduct} we introduce the product $\smile$ on
$\rous(Y)$ for any stably complex $T$-manifold $Y$.
Its definition will be of no surprise to anyone who has computed with
Chen-Ruan cohomology, though it does not seem to have been formalized
before in terms of equivariant cohomology; it is set up to make it
easy to show that
$$ \rous(\locfree)\iso H_{CR}^*(\locfree/T)
$$
for any space $\locfree$ with a locally free $T$-action.

We then prove the essential fact that in the robustly
equivariantly injective case, the product $\star$ can be used to
compute the product $\smile$. The virtue of $\star$ is that it is
easy to compute with; for example, it is essentially automatic to
show that it is graded and associative. On the other hand, we need
the $\smile$ product
% (even in the robustly equivariantly injective case)
to show that we  indeed have a well-defined product on $\rous(Y)$,
rather than merely on $\rousT$.
In addition, $\smile$ is better behaved from a
functorial point of view, such as when restricting to a level set of the
moment map of a proper Hamiltonian $T$-space.

%, or proving equivariant versions of our main theorems.

%\begin{theorem}
%Let $Y$ be a Hamiltonian $T$-space with a proper moment map bounded from below. There is a natural surjection of %rings
%$$
%\rous(Y)\longrightarrow H_{CR}(Y/\!/T(\mu))
%$$ for any regular value $\mu$ of the moment map.
%\end{theorem}
%We also describe the kernel of this surjection.

Suppose that $Y$ is a proper Hamiltonian $T$-space, and $\locfree$ is
the zero level set $\Phi^{-1}(\mu)$.
One of our main theorems states that there is a surjection of graded
(in the first coordinate) rings
\begin{equation}\label{eq:surjection}
  \rous(Y) \onto H_{CR}^*(Y/\!/T(\mu))
\end{equation}
arising from a natural restriction map $\rous(Y) \longrightarrow
\rous(\locfree)$. This follows from the work of Kirwan
\cite{K:quotients} and the fact that
$\rous(\locfree)\iso H_{CR}^*(\locfree/T)$.
Furthermore, the ring $\rous(Y)$ is {\em easy} to compute: the
required data is readily available from the symplectic point of view.
The kernel of the map (\ref{eq:surjection}) may be computed using
techniques introduced by Tolman and Weitsman \cite{TW:symplecticquotients}
and refined by the first author in \cite{Go:effective}.
Essentially, our definitions and theorems are generalizations to
Chen-Ruan cohomology of similar ones about
the maps
$$H^*_T(Y) \into H^*_T(Y^T), \qquad H^*_T(Y) \onto H^*(Y/\!/T) $$
familiar in Hamiltonian geometry.
An easy observation is that in order to compute $H_{CR}^*(Y/\!/T(\mu))$ for any particular value $\mu$, one only needs a much smaller ring. Let $\Gamma_\mu\subset T$ be the subgroup generated by all finite stabilizers occurring in the $T$ action on $\Phi^{-1}(\mu)$. Then
$$
NH_T^{*,\Gamma_\mu}(Y) := \bigoplus_{g\in \Gamma_\mu}  H_T^*(Y^g)
$$ is a subring of $\rous(Y)$ that also surjects onto $H_{CR}^*(Y/\!/T(\mu))$. In particular, if $\Gamma$ is the subgroup generated by all finite stabilizers occurring in the $T$ action on $Y$, then $NH_T^{*,\Gamma}(Y)$ surjects onto $H_{CR}^*(Y/\!/T(\mu))$ for every regular value $\mu$.

We consider the whole ring $\rous(Y)$ rather than just this subring because we think it is interesting in its own right; it also lends elegance to proofs and statements of results. While $\Gamma$ is easily computed from $T$ acting on $Y$, it is an unnecessary computational step in order to state the surjectivity result (\ref{eq:surjection}).  In addition, there is no natural map on inertial cohomology given a homomorphism of groups $\Gamma_1\rightarrow \Gamma_2$. For example, Proposition~\ref{prop:transversehom} would not hold if the finite stabilizers occurring on $X$ were different from those occurring on $Y$. By allowing the group over which we take the direct sum to be as large as possible, we are able to obtain some results concerning the functoriality of $\rous$.

 The paper is organized as follows. In
Section~\ref{sec:starproduct} we define $\rous(Y)$ as an
$H_T^*(pt)$-module, the ``restriction'' map from $\rous(Y)$ to
$H^*_T(Y^T) \tensor_{\Z} \Z[T]$, and the product $\star$ on $H^*_T(Y^T)
\tensor_{\Z} \Z[T]$. We show that $\star$ is an associative product,
and graded. At this stage it is unclear that the image of the
restriction map is closed under $\star$ (and this will wait until
Section \ref{sec:smileproduct}), but if one accepts this as a black
box one can already begin computing examples. When $Y$ is robustly
equivariantly injective,
%R: I DON'T THINK THE READER GETS ANYTHING BUT CONFUSION FROM THIS COMMA-ED REMARK.
%ALSO, I MADE THIS SENTENCE PART OF THE PRECEDING PARAGRAPH.
%so the restriction map
%is injective, and assuming the $\star$-closedness black box,
the product $\star$ induces a product on $\rous(Y)$.

In Section~\ref{sec:smileproduct} we define the $\smile$ product on $\rous(Y)$
for any stably complex manifold with a smooth $T$-action.
This definition makes the grading $\diamond$ over elements of $T$ obvious
while obscuring the associativity and the grading by real numbers.

Our main theorem in this section, Theorem \ref{th:productsthesame},
is that the restriction map (\ref{eq:introinjective}) is a
%R: I ADDED AN EQUATION NUMBER TO THE RESTRICTION MAP
%I CUT THE NOT NECESSARILY ASSOCIATIVE REMARK BECAUSE IT IS CONFUSING TO A READER. THE MAP IS ASSOCIATIVE,
%JUST THAT WE SHOW IT IS A RING HOMOMORPHISM BEFORE SHOWING ASSOCIATIVE. WE MAKE THIS POINT IN THE TEXT
%ANYWAY AND IT SEEMS MISLEADING, AS IF THIS WERE AN IMPORTANT POINT, IN THE INTRO.
%(not necessarily associative)
ring homomorphism from $\smile$ to $\star$.
In particular, the image is a subring, and when $Y$ is robustly
equivariantly injective, the $\star$ product can be used as a simple
means of computing the ring $\rous(Y)$. For example, the associativity and
gradedness of $\star$ prove the same properties of $\smile$. (In
fact $\smile$ has these properties even when $Y$ is not robustly
equivariantly injective.)

In Section~\ref{sec:orbifoldcohomology}
we prove that the inertial cohomology (with the $\smile$ product)
of a space $\locfree$ with a locally free $T$-action
is isomorphic to the Chen-Ruan cohomology ring of the quotient orbifold;
this is essentially a definition chase and was our motivation for $\smile$.
We also show that, for a stably almost complex
manifold carrying a $T$-invariant function, the inclusion of a regular level
set induces a well-defined map in inertial
cohomology. As a corollary we obtain surjectivity from the inertial
cohomology ring of a proper Hamiltonian $T$-space
%(equipped with a moment map with one component proper and bounded below),
to the Chen-Ruan cohomology ring of the symplectic reduction. This
connection is elaborated upon in Section~\ref{sec:surjectivity}.
%There is also an equivariant version of this surjectivity, which relies on
%the functorial properties of $\rous(Y)$, described in
%Section~\ref{sec:functoriality}.
Finally, we spend significant effort in making these computations tenable.
In Section~\ref{sec:multiplication} we give yet another description
of the product, and in Sections~\ref{sec:singularities} and
\ref{sec:toricvarieties} we explore two important sets of examples,
namely weight varieties (symplectic reductions of coadjoint orbits)
and symplectic toric orbifolds.

Since completing this work (announced in \cite{GHK:ober}), we received the
preprint \cite{CH:deRham}, which also uses equivariant cohomology to 
study Kirwan surjectivity for abelian symplectic quotient orbifolds
(though it doesn't address Kirwan injectivity).
Their introduction of ``twist factors'' (Definition 3.1 in \cite{CH:deRham}) 
into the de Rham models of ordinary and equivariant cohomology parallels 
closely our modified homomorphism in Section \ref{sec:multiplication}.
A crucial ingredient in both our approach and theirs
is a simplification in the description of the obstruction bundle
(presented here in Definition \ref{de:obbundle});
we got it from \cite[Proposition 6.3]{BCS:toricvarieties},
whereas they rederive it in \cite[Proposition 3.4]{CH:deRham},
just before Equation (3.7).
Then the work in our central Theorem \ref{th:productsthesame} 
takes place in their Propositions 3.2 and 3.4. 
The fact that their twist factors are cohomology classes rather than differential forms 
(in their words, ``to avoid unnecessary nonsense of choosing forms'') 
makes it appear more natural to work in cohomology from the start, 
as we do, which also allows us to work over $\Z$ rather than $\R$.

%\todo{The paragraph above is what should be elaborated on in detail -
%  the editor wants something serious it seems! I think we should keep
%  such an elaborated version in the intro, but then if necessary put
%  in a comment here or there in the body of the work at points where
%  the work or definitions are similar.-R}

%%%%%%%%%%%%%%%%%%%%%%%%%%%%%%%%%%%%%%%%%%%%%%
%We begin by defining a vector space $\gpiece(Y)$ for any $g\in T$. As a vector space, and indeed as a module over $H_T^*(pt)$, let $\gpiece(Y):= H_T^*(Y^g)$.  The grading on $\gpiece(Y)$ is shifted, however. Note that $g$ acts on the normal bundle to $Y^g$ in $Y$. For any connected component $Z$ of $Y^g$, this action splits $\nu Z$ into a sum of isotypic components $\oplus_{\lambda}I_\lambda$, where $g$ acts on $I_\lambda$ with  eigenvalue $e^{2\pi ia_\lambda}$, where $a_\lambda\in [0,1)$. Then  the {\em age} of $g$ on $Z$ is given by
%$$
%\age(Z,g):= \sum_\lambda \dim(I_\lambda)a_\lambda.
%$$
%Let $\alpha\in H_T^j(Y^g)$ be such that $\alpha|_Z\neq 0$ and $\alpha|_{Z'}=0$ for any connected components $Z'\neq %Z$ in $Y^g$. Then by definition,
%$$
%\deg(\alpha) = j + 2\age(Z,g).
%$$

%The grading on $\rous(Y)$ is inherited from that on $\gpiece(Y)$ for each $g$.
%All terms in $\rous(Y)$ that are not homogeneous both in the first and second gradings may be written as a sum of classes that are homogeneous in each grading, and which may be in turn written as a sum of classes each of which is supported on exactly one connected component of some piece $Y^g$ of $Y$.

{\bf Acknowledgements.} We are very thankful to Tom Graber for
demystifying Chen-Ruan cohomology for us, and to Daniel Biss for
pointing out that we define a Ring of Unusual Size \cite{PB}. We also
thank Megumi Harada for a careful reading of earlier drafts.

%%% Local Variables: 
%%% mode: latex
%%% TeX-master: t
%%% End: 

%% 2. Preorbifold cohomology for Hamiltonian T-spaces -- Rebecca
%%    hamTsp.tex
%%        * Define via equivariant injectivity

\section{Inertial cohomology and the $\star$ product}
\label{sec:starproduct}
Let $Y$ be a stably almost complex $T$-space.
% robustly equivariantly injective $T$ space.
For any $g\in T$, let $Y^g :=\{y\in Y\ |\ g\cdot y = y\}$ denote the set of
points fixed by $g$. Denote by $\stab(y)$ the stabilizer of $y$ in $T$.
% We begin by defining the {\dfn $g$-piece} of the inertial
% cohomology $\rous(Y)$, in which the grading in $\diamond$ is implicit.
% AK: this seems an unnecessary distraction, as I found it used in exactly
% one place, and could be safely defined there if need be.
Since $T$ is abelian, each $Y^g$ is itself naturally a $T$-space.

\begin{definition}
  The {\dfn inertial cohomology} of the space $Y$ is given as an
  $H_T^*(pt)$-module by
  \begin{align*}
    \rous(Y) :&= \bigoplus_{g\in T} H_T^*(Y^g)
  \end{align*}
  where the sum indicates the $\diamond$ grading,
  i.e. $\gpiece(Y) = H_T^*(Y^g)$.
\end{definition}

% For a class $b\in\rous(Y)$, let $b_g$ denote the component of $b$
% living in $\gpiece(Y)$.
% We say that $b\in \gpiece(Y)$ or $b$ is
% {\dfn supported} on $Y^g$ if $b_h=0$ for $h\neq g$.
% AK: this is a goofy definition, since Y^g may equal Y^h.
% Moreover, we use b_1..b_n below for something else.
Notice that if $Y$ is compact and $T$ acts on $Y$ locally freely, i.e.
if $\stab(y)$ is finite for all $y\in Y$, then there are only finitely
many nonzero terms in the sum above.
% defining sum $\rous(Y)$ is a
% finite sum of elements in $T$ appearing as stabilizers.
At the other extreme, if $T$ has fixed points $Y^T$ on $Y$, then every
summand is nonzero.

Neither the ring structure nor the first grading ${}^*$
are the ones induced from each $H_T^*(Y^g)$.
They will be defined in Section \ref{sec:smileproduct},
and will depend on $Y$'s stably almost complex structure
(which the definition above does not).

For $L$ a $1$-complex-dimensional representation of $T$ with weight $\lambda$,
and $g\in T$, the eigenvalue of $g$ acting on $L$ is
$\exp({2\pi i\,\Year_{\lambda}(g)})$ where
$\Year_{\lambda}(g) \in [0,1)$ is the {\dfn logweight} of $g$ on $L$.

Since $T$ preserves the stably almost complex structure on $Y$,
any component $F$ of the fixed point set $Y^g$ is also stably almost complex,
and the normal bundle $\nu F$ to $F$ in $Y$ is
an {\em actual} complex vector bundle. The torus $T$ acts on $\nu F$,
and splits it into isotypic components
$$
\nu F =\bigoplus_\lambda I_\lambda
$$
where the sum is over weights $\lambda\in \widehat{T}$.
We denote the logweight of $g$ on $I_\lambda$
(restricted to any point $y\in F$) by $\Year^F_\lambda(g)$.

For each $g\in T$, there is an inclusion $Y^T \into Y^g$, inducing
a map backwards in $T$-equivariant cohomology. (For most $g$,
this inclusion is equality.) Put together, these give a {\dfn restriction} map
$$ i_{NH}^*:
\rous(Y) \longrightarrow \bigoplus_{g\in T} H^*_T(Y^T) \iso \rousT
$$
where the isomorphism is as $H^*_T(pt)$-modules. If there are no
fixed points, this map is zero. The most interesting case is when
$Y$ is robustly equivariantly injective, meaning that this map is
injective.

% AK moved a bunch of stuff to leftovers here

For $b\in \rousT$, $F$ a component of $Y^T$, and $g\in T$,
let $b|_{F,g} \in H^*_T(F)$ denote the {\dfn component} of $b$
in the $g$ summand, restricted to the fixed-point component $F$.
Only finitely many of these components can be nonzero,
and $b$ can be reconstructed from them as
$$ b = \sum_{F,g} (b|_{F,g}) \tensor {1g} \quad \in \rousT.$$
where $b|_{F,g} \in H^*_T(F)$ lives in $H^*_T(Y^T)$
via the decomposition $H^*_T(Y^T) = \bigoplus_{F_i} H^*_T(F_i)$.

\begin{definition}\label{def:starproduct}
  Let $Y$ be a stably almost complex $T$-space.

  Let $b_1,b_2 \in \rousT$. The product
  $b_1\star b_2$ is defined componentwise:
  \begin{equation}\label{eq:starproduct}
    (b_1\star b_2)|_{F,g} :=
    \sum_{(g_1,g_2) : \atop g_1 g_2 = g}
    (b_1|_{F,g_1}) (b_2|_{F,g_2})
    \prod_{I_{\lambda}\subset \nu F}
    e(I_\lambda)^{\Year^F_\lambda(g_1)+\Year^F_\lambda(g_2)
      -\Year^F_\lambda(g_1 g_2)} %\in  NH_T^{*,g_1\cdots g_n}(F)
  \end{equation}
  where $e(I_\lambda)\in H_T^*(F)$ is the equivariant Euler class of
  $I_\lambda$.

  More generally, we have an $n$-ary product:
  \begin{equation}\label{eq:starproductn}
    (b_1\star\cdots \star b_n)|_{F,g} :=
    \sum_{(g_1,\ldots,g_n) : \atop \prod g_i = g}
    \prod_i (b_i|_{F,g_i})
    \prod_{I_{\lambda}\subset \nu F}
    e(I_\lambda)^{\Year^F_\lambda(g_1)+\ldots+\Year^F_\lambda(g_n)
      -\Year^F_\lambda(g_1g_2\cdots g_n)}. %\in  NH_T^{*,g_1\cdots g_n}(F)
  \end{equation}
\end{definition}
Note that the exponent in (\ref{eq:starproductn}) is an integer from 0
to $n-1$, and in fact is the greatest integer
$\big\lfloor \Year^F_\lambda(g_1)+\ldots+\Year^F_\lambda(g_n) \big\rfloor$.
% We prove in Theorem~\ref{th:associative} that $\star$ is an
% associative product.
Since the sum is over $\{(g_1,\ldots,g_n) : \prod g_i = g\}$,
this product plainly respects the $T$-grading from the second
factor of $\rousT$.

\begin{theorem}\label{th:associative}
  The $2$-fold product $\star$ is associative, making $\rousT$ into a ring.
\end{theorem}

\begin{proof}
  Let $b_1,b_2,b_3 \in \rousT$.
  We relate the $2$-fold product to the $3$-fold:
  $$ (b_1 \star b_2) \star b_3 = b_1 \star b_2 \star b_3
  = b_1 \star (b_2 \star b_3).$$
  (More generally, the $n$-fold product can be built from the $2$-fold
  with any parenthesization.)

  There are two multiplicative contributions to a $\star$-product:
  the components $b_i|_{F,g}$, and the equivariant Euler classes
  of the $I_\lambda$.
  For a given $F$ and triple $g_1,g_2,g_3$ of group elements, the components
  give the same contribution
  $$
    (b_1|_{F,g_1}\cdot b_2|_{F,g_2})\cdot b_3|_{F,g_3}
  =  b_1|_{F,g_1}\cdot b_2|_{F,g_2} \cdot b_3|_{F,g_3}
  =  b_1|_{F,g_1}\cdot(b_2|_{F,g_2} \cdot b_3|_{F,g_3}).
  $$
  To see that the exponents match on the equivariant Euler class of
  $I_\lambda$, we need to check
  \begin{align*}
&\big[\Year^F_\lambda(g_1)+\Year^F_\lambda(g_2)-\Year^F_\lambda(g_1g_2))\big] +
 \big[\Year^F_\lambda(g_1g_2)+\Year^F_\lambda(g_3)-\Year^F_\lambda((g_1g_2)g_3)\big]\\
     =\ &    \Year^F_\lambda(g_1)+\Year^F_\lambda(g_2)+\Year^F_\lambda(g_3)-\Year^F_\lambda(g_1g_2g_3)  \\
     =\ &
    \big[\Year^F_\lambda(g_2)+\Year^F_\lambda(g_3)-\Year^F_\lambda(g_2g_3)\big]+\big[\Year^F_\lambda(g_1)+
    \Year^F_\lambda(g_2g_3)-\Year^F_\lambda(g_1(g_2g_3))\big]
  \end{align*}
  which is plain.
\end{proof}

\begin{remark}\label{rem:idprod}
  If $b_1 \in H^*_T(Y^T)\tensor {\bf 1}, b_2 \in H^*_T(Y^T)\tensor g$,
  where ${\bf 1}$ is the identity element of $T$,
  then $b_1 \star b_2 \in H^*_T(Y^T)\tensor g$, and
  $$ (b_1 \star b_2)|_F = (b_1|_{F,{\bf 1}}) (b_2|_{F,g})
    \prod_{I_{\lambda}\subset \nu F}
    e(I_\lambda)^{\Year^F_\lambda({\bf 1})+\Year^F_\lambda(g)
      -\Year^F_\lambda(g)} %\in  NH_T^{*,g_1\cdots g_n}(F)
    = (b_1|_{F,{\bf 1}}) (b_2|_{F,g})
    \prod_{I_{\lambda}\subset \nu F}
    e(I_\lambda)^0
    = (b_1|_{F,{\bf 1}}) (b_2|_{F,g})
  $$
  is the ordinary multiplication.
\end{remark}

What is not clear at this point is that
the image of the restriction map from $\rous(Y)$ into $\rousT$
is closed under $\star$.
This will follow from Theorem \ref{th:productsthesame} in the next section.

We now turn to the grading in $\rousT$, which is built
from logweights. In particular, it is not in general graded by integers,
nor even by rationals (unlike Chen-Ruan cohomology).

Let $y\in Y^g$, and decompose $T_y Y=\oplus_j L_j$ under the $g$ action.
The sum of the logweights of $g$ on each of these lines is termed the
{\dfn age} of $g$ at $y$ (see \cite{R:age}). Since this number depends
only on the connected component $Z$ of $y$ in $Y^g$, let
$$
\age(Z,g) = \sum_j \Year_{\lambda_j}(g).
$$
It is in general a real number. If $g$ is of finite order $n$, then
$\age(Z,g) \in \frac{1}{n}\Z$, but our $g$ are not in general of finite order. See Remark~\ref{re:rationalage}.

Let $b\in \rousT$ such that $b|_{Z,g}$ is zero except for one $g \in T$
and one component $Z$ of $Y^g$. Further assume that
$b|_{Z,g} \in H^*_T(Z)$ is homogeneous.
Then we assign
$$ \deg(b)=\deg(b|_{Z,g}) + 2\ \age(Z,g).$$
(The $2$ is the usual conversion factor $\dim_\R \C$; if one works with
Chow rings rather than cohomology one doesn't include it.)
%Note that the grading is {\em real} rather than integral.

\begin{theorem}\label{th:graded}
  This definition of degree makes $(\rousT,\ast)$ into a graded algebra.
\end{theorem}

\begin{proof}
  It is enough to check when $b_1,b_2$ each have only one nonvanishing
  component $b_1|_{F_1,g_1}$, $b_2|_{F_2,g_2}$, hence degrees
  $$ \deg(b_i)=\deg(b_i|_{F_i,g_i}) + 2\ \age(F_i,g_i).
  $$
  If $F_1 \neq F_2$, then $b_1 \star b_2 = 0$ and there is nothing to prove.
  Otherwise $b_1 \star b_2$ has only one nonvanishing component,
  $$ (b_1 \star b_2)|_{F,g_1 g_2} = b_1|_{F,g_1} b_2|_{F,g_2}
  \prod_{I_{\lambda}\subset \nu F}
  e(I_\lambda)^{\Year^F_\lambda(g_1)+\Year^F_\lambda(g_2)
    -\Year^F_\lambda(g_1g_2)} $$
  where $F = F_1 = F_2$. The degree we wish to assign this is
  \begin{align*}
   &\deg (b_1 \star b_2)  \qquad\qquad\qquad\qquad\qquad \\
   &= \bigg( \deg b_1|_{F,g_1} + \deg b_2|_{F,g_2}
  + \sum_{I_{\lambda}\subset \nu F}
  \big({\Year^F_\lambda(g_1)+\Year^F_\lambda(g_2)
    -\Year^F_\lambda(g_1g_2)}\big) \   (2 \dim_\C I_\lambda) \bigg)
  + 2\ \age(F,g_1 g_2).
  \end{align*}
  Canceling the $b_i|_{F,g_i}$ contributions and the factor of $2$,
  our remaining task is therefore to show that
  $$ \age(F,g_1) + \age(F,g_2) =
  \age(F,g_1 g_2) +
  \sum_{I_{\lambda}\subset \nu F}
  \big({\Year^F_\lambda(g_1)+\Year^F_\lambda(g_2)
    -\Year^F_\lambda(g_1g_2)}\big) \   \dim_\C I_\lambda .
  $$
  This follows from three applications of the formula
  $$ \age(F,g) = \sum_{I_{\lambda}\subset \nu F}
  \Year^F_\lambda(g) \   \dim_\C I_\lambda$$
  which is just a resummation of the definition.
% AK: original proof moved to leftovers
\end{proof}

Our point of view is that the grading is not of fundamental importance --
it happens to be preserved by the multiplication, so we record it as
an extra tool for studying this ring.
In Section \ref{sec:multiplication} we will see what is perhaps the best
motivation for this grading.

% \todo{We could do a simple application here, confirming that the image
%   of restriction is indeed a subring, e.g. $S^1$ acting on $\P^2$ with
%   weights $0,1,2$. This would emphasize that calculation is already
%   possible. If we do, we have to change the first paragraph in the
%   next section, which presently says that we wait until the end to give
%   examples.}

%% 3. Definition of preorbifold cohomology -- Rebecca
%%    defn.tex
%%        * The ROUS def'n
%%        * The nice finite subring (the \Gamma-subring)
%%        * Thm: The ROUS is a graded associative ring
%%        * Thm: This is the same as the previous defn

\section{Inertial cohomology and the $\smile$ product}
\label{sec:smileproduct}
The definition of the $\star$ product on $\rousT$ renders the ring
$\rous(Y)$ straightforward to compute when $Y$ is robustly
equivariantly injective, as will be shown in the examples at the end
of this paper. It also has the advantage that associativity is easy
to prove. However, its limits are easy to see: if there are no fixed
points, for example, $\rousT$ is zero; we have no proof yet that
$i_{NH}^*(\rous(Y))$ is a subring of $\rousT$; even assuming that,
when $Y$ is not robustly equivariantly injective, the $\star$
product doesn't let us define a product on the source $\rous(Y)$,
only on its image inside $\rousT$.
% On the other hand, suppose $Y$ is proper Hamiltonian, which
% ensures that $\star$ has a rich structure. Then the level set of the
% moment map is not itself Hamiltonian (or equivariantly injective).
For these reasons, we present in this section a product $\smile$
directly on $\rous(Y)$ for any stably almost complex manifold $Y$,
and show that $i_{NH}^*$ is a ring homomorphism.
The $\smile$ product has its roots in the
original paper by Chen and Ruan \cite{CR:orbH}, but
is defined using the global group action and the language of
equivariant cohomology.

\newcommand\tildeY{\widetilde Y}
\newcommand\barE{\overline{e}}

Let $Y$ be a stably almost complex manifold
with a smooth $T$ action respecting the stably almost complex
structure. Note that this implies that each normal bundle
$\nu(Y^{g_1,g_2}\subset Y^{g_1})$ is a complex vector bundle
(not just stably so) over
$Y^{g_1,g_2}$ for every choice of $g_1,g_2\in T$, where
$Y^{g_1,g_2}=(Y^{g_1})^{g_2}=(Y^{g_2})^{g_1}$.
The definition of the product $\smile$ requires the introduction of a
new space, and a vector bundle over each of its connected components.
Let
$$
\tildeY:= \coprod_{g_1,g_2\in T} Y^{g_1,g_2}.
$$
For any connected component $Z$ of $Y^{g_1,g_2}$, the group $\langle
g_1,g_2\rangle$ generated by $g_1$ and $g_2$ acts on the
complex vector bundle $\nu Z$, the normal bundle to $Z$ in $Y$, fixing
$Z$ itself. Thus as a representation of $\langle g_1,g_2\rangle$, $\nu
Z$ breaks up into isotypic components
$$
\nu Z = \bigoplus_{\lambda\in \widehat{\langle
g_1,g_2\rangle}} I_\lambda
$$
where $I_\lambda$ is the bundle over $Z$ on which $\langle
g_1,g_2\rangle$ acts with representation given by $\lambda$.
%We define
%the {\em time of $I_\lambda$} to be the sum
%$$
%\Time_\lambda(g,h) =
%\Year_\lambda(g)+\Year_\lambda(h)+\Year_\lambda((gh)^{-1}).
%$$ We will alternatively denote this by $\Time(I_\lambda)$, provided the omitted
%information is clear. Note that $\Time(I_\lambda)$ is 0, 1, or 2.

\begin{definition}\label{de:obbundle}
  For each connected component $Z$ of $Y^{g_1,g_2}$ in $\tildeY$,
  let $E|_Z$ be the vector bundle over $Z$ given by
  $$
  E|_Z =
  \bigoplus_{\Year_\lambda(g_1)+\Year_\lambda(g_2)+\Year_\lambda(g_3)=2}
  I_{\lambda},
  $$ where $g_3:= (g_1g_2)^{-1}$.
                                %\Time_\lambda(g,h)=2
  The {\dfn obstruction bundle} $E$ is given by the union of $E|_Z$ over
  all connected components $Z$ in $\tildeY$.
\end{definition}

This sum of three logweights is reminiscent of, but {\em not} like, age:
one calculates the age of a group element by summing over lines, whereas
this functional is calculated for a line by summing over group elements.

Note that the dimensions of the fibers of $E$ may differ on different
connected components.

\begin{remark}
  Each component $Z$ is $T$-invariant and hence $E|_Z\to Z$ is a
  $T$-equivariant bundle. Thus there is a well-defined (inhomogeneous)
  equivariant Euler class $\varepsilon$ of $E$: for every component
  $Z$, let $\varepsilon$ restricted to $Z$ be the equivariant Euler
  class of $E|_Z$. The class $\varepsilon$ is also called the
  {\dfn virtual fundamental class} of $\tildeY$.
\end{remark}

Consider the three inclusion maps given by
\begin{align*}
e_1: &Y^{g_1,g_2}\hookrightarrow Y^{g_1}\\
e_2: &Y^{g_1,g_2}\hookrightarrow Y^{g_2},\ \mbox{ and}\\
\barE_3: &Y^{g_1,g_2}\hookrightarrow Y^{g_1g_2},
\end{align*}
(The notation will be explained in Section \ref{sec:orbifoldcohomology}.)
The maps $e_1,e_2,\barE_3$ clearly extend to maps on
$\tildeY$. They therefore induce the pullbacks
$$
e_1^*, e_2^*: \rous(Y)\To \bigoplus_{g_1,g_2\in T}H_T^*(Y^{g_1,g_2})
$$
and the pushforward map
$$
(\barE_3)_*: \bigoplus_{g_1,g_2\in T}
H_T^*(Y^{g_1,g_2})\To \rous(Y).
$$
\begin{definition}
For $b_1,b_2\in \rous(Y)$, define
$$
b_1\smile b_2 := (\barE_3)_*\ \big(e_1^*(b_1)\cdot
e_2^*(b_2)\cdot  \varepsilon\big),
$$
where $\varepsilon$ is the virtual fundamental class of $\tildeY$,
and the product occurring in the right hand side is that in the
equivariant cohomology of each piece $Y^{g_1,g_2}$ of $\tildeY$.
\end{definition}

\begin{remark}
  While $Y^{g_1,g_2}$ and $Y^{g_1 g_2}$ may be noncompact, $\barE_3$
  is still a closed embedding, so the pushforward is well-defined.
\end{remark}

\begin{remark}\label{re:whereproductlives}
If $b_1\in NH_T^{*,g_1}(Y)$ and $b_2\in NH_T^{*,g_2}(Y)$, then
$e_1^*(b_1)$, $e_2^*(b_2)$, and their product
$e_1^*(b_1)\cdot e_2^*(b_2)$ live in $H_T^*(Y^{g_1,g_2})$.
After multiplying by $\varepsilon$, the pushforward map
$(\barE_3)_*$ sends this class to $H_T^*(Y^{g_1g_2})$, which
implies that $b_1\smile b_2\in NH_T^{*,g_1g_2}(Y)$.
\end{remark}

We also define a real-valued grading on $\rous(Y)$. Recall that
$$ NH^{*,g}_T(Y) = H^*_T(Y^g) = \bigoplus_Z H^*_T(Z) $$
where $Z$ varies over the connected components of $Y^g$.
We shift the degree on the $H^*_T(Z)$ summand by
(twice) the age of $g$ on any tangent space $T_z Y, z\in Z$.
In particular we usually
do {\em not} shift all of $H^*_T(Y^g)$ by the same amount.

At this point we haven't shown that $\smile$ is associative, so the
main theorem of this section has to be phrased in terms of
``not-necessarily-associative rings.''

\begin{theorem}\label{th:productsthesame}
  Let $Y$ be a stably almost complex manifold with $T$ action,
  and let $b_1,b_2 \in \rous(Y)$. Let $i_{NH}^*:\rous(Y)\rightarrow \rousT $
  be the restriction map.
  Then
  $$ i_{NH}^*(b_1 \smile b_2) = i_{NH}^*(b_1) \star i_{NH}^*(b_2), $$
  i.e. $i_{NH}^*$ is a homomorphism of not-necessarily-associative rings.
  Moreover, it preserves the real-valued grading.
\end{theorem}

\begin{proof}
  Suppose that $b_1\in NH_T^{*,g_1}(Y)$ and $b_2\in NH_T^{*,g_2}(Y)$
  (it is enough to prove it for this case).
  Our goal is to prove
  $$ i_{NH}^*(b_1 \smile b_2)|_{F,g}
  = \big(i_{NH}^*(b_1) \star i_{NH}^*(b_2)\big)|_{F,g} $$
  for every $(F,g)$ component of either side. It is easy to see that
  both sides vanish unless $g = g_1 g_2$.

  Let $F$ be a connected component of the fixed point set $Y^T$,
  which therefore includes $T$-equivariantly into $Y^{g_1,g_2}$,
  which in turn includes into each of $Y^{g_1}$, $Y^{g_2}$, $Y^{g_1 g_2}$.
  Call this first map $i_F : F \to Y^{g_1,g_2}$, the others already
  having names $e_1,e_2,\barE_3$.

  We start with the left side. Let $Z$ denote the connected component of
  $Y^{g_1,g_2}$ containing $F$, let $\epsilon_Z$ denote the Euler class
  of the obstruction bundle on $Z$, and let $f_Z$ denote the equivariant
  Euler class of $Z$'s normal bundle inside $Y^{g_1 g_2}$. Then
  \begin{align*}
    i_{NH}^*(b_1 \smile b_2)|_{F,g_1 g_2}
    &=& (\barE_3 \circ i_F)^*
    \bigg((\barE_3)_* (e_1^*(b_1) \cdot e_2^*(b_2) \cdot \epsilon_Z)\bigg)
    &\hfill& \hbox{by the definition of $\smile$} \\
    &=&  i_F^* \barE_3^*
    \bigg((\barE_3)_* (e_1^*(b_1) \cdot e_2^*(b_2) \cdot \epsilon_Z)\bigg)&&\\
    &=&  i_F^* \bigg(e_1^*(b_1) \cdot e_2^*(b_2) \cdot \epsilon_Z
    \cdot \barE_3^* ((\barE_3)_* 1)\bigg)
    &\hfill& \hbox{by the pull-push formula} \\
    &=&  i_F^* e_1^*(b_1) \cdot i_F^* e_2^*(b_2)
    \cdot i_F^* \epsilon_Z \cdot i_F^* f_Z
    &\hfill& \hbox{by the definition of Euler class.} \\
  \end{align*}
  Now we compare to the right side.
  \begin{align*}
    \big(i_{NH}^*(b_1) \star i_{NH}^*(b_2)\big)|_{F,g_1 g_2}
    &=& i_{NH}^*(b_1)|_{F,g_1} \cdot i_{NH}^*(b_2)|_{F,g_2} \cdot
    \prod_{I_{\lambda}\subset \nu F}
    e(I_\lambda)^{\Year^F_\lambda(g_1)+\Year^F_\lambda(g_2)
      -\Year^F_\lambda(g_1g_2)} \\
    &=& i_F^* e_1^*(b_1) \cdot i_F^* e_2^*(b_2) \cdot
    \prod_{I_{\lambda}\subset \nu F}
    e(I_\lambda)^{\Year^F_\lambda(g_1)+\Year^F_\lambda(g_2)
      -\Year^F_\lambda(g_1g_2)}
  \end{align*}
  Our goal is thus to show that
  $$ i_F^* \epsilon_Z \cdot i_F^* f_Z =
  \prod_{I_{\lambda}\subset \nu F}
  e(I_\lambda)^{\Year^F_\lambda(g_1)+\Year^F_\lambda(g_2)
    -\Year^F_\lambda(g_1g_2)}.
  $$
  In fact we will show by case analysis that for each
  $\lambda \in \widehat{\langle g_1,g_2\rangle}$, the equivariant Euler
  class of the bundle $I_\lambda$ over $F$ shows up to the same power
  on the left and the right side. Let
  \begin{eqnarray*}
    \varepsilon(\lambda) &=&\begin{cases}
      1 &\mbox{if
        $\Year^F_\lambda(g_1)+\Year^F_\lambda(g_2)+\Year^F_\lambda(g_3)=2$,
        where $g_3=(g_1g_2)^{-1}$}\\
      0 &\mbox{otherwise, and}
    \end{cases}\\
    \mathnormal{f}(\lambda)&=&\begin{cases}
      1
      &\mbox{if $I_\lambda\subset \nu(Y^{g_1,g_2} \subset Y^{g_1g_2})$}\\
      0 &\mbox{otherwise.}
    \end{cases}
  \end{eqnarray*}
  So $i_F^* \epsilon_Z = \prod_\lambda e(I_\lambda)^{\varepsilon(\lambda)}$ and
  $i_F^* f_Z = \prod_\lambda e(I_\lambda)^{f(\lambda)}$, and it
  remains to check that
  $$ \varepsilon(\lambda) + f(\lambda) =
  \Year^F_\lambda(g_1)+\Year^F_\lambda(g_2) - \Year^F_\lambda(g_1 g_2).
  $$
  The main principle in the following case analysis is that
  $\Year^F_\lambda(g_1)+\Year^F_\lambda(g_2) - \Year^F_\lambda(g_1 g_2)$
  is either $0$ or $1$, not some arbitrary real number, and likewise
  $\Year^F_\lambda(g_1)+\Year^F_\lambda(g_2) + \Year^F_\lambda(g_3)$
  is either $0$, $1$, or $2$.

  Assume first that $\Year^F_\lambda(g_1 g_2) = 0$,
  meaning $g_1 g_2$ acts trivially on $I_\lambda$.
  Then $\Year^F_\lambda(g_3) = 0$, hence $\varepsilon(\lambda) = 0$;
  also $I_\lambda \leq TY^{g_1 g_2}$.
  So the equation we seek is
  $f(\lambda) = \Year^F_\lambda(g_1)+\Year^F_\lambda(g_2)$, where both
  sides are either $0$ or $1$. Now
  $$ f(\lambda)=0
  \qquad \hbox{iff}\qquad
  \hbox{$g_1,g_2$ each act trivially on $I_\lambda$}
  \qquad \hbox{iff}\qquad
  \Year^F_\lambda(g_1) = \Year^F_\lambda(g_2) = 0. \qquad \hbox{QED.}$$

  On the other hand, assume that $\Year^F_\lambda(g_1 g_2) \neq 0$.
  Then $I_\lambda \not\leq TY^{g_1 g_2}$, hence $f(\lambda) = 0$.
  So the equation we want now is
  $\varepsilon(\lambda)
  = \Year^F_\lambda(g_1)+\Year^F_\lambda(g_2) - \Year^F_\lambda(g_1 g_2)$.
  The right side is $1$ iff
  $\Year^F_\lambda(g_1)+\Year^F_\lambda(g_2) > 1$ iff
  $\Year^F_\lambda(g_1)+\Year^F_\lambda(g_2) + \Year^F_\lambda(g_3) = 2$
  iff $\varepsilon(\lambda) = 1$. QED.

  It is trivial to check that the grading is respected,
  essentially because both gradings are defined using ages.
  (In particular, the proof does not require splitting into cases.)
\end{proof}

If $Y$ is robustly equivariantly injective, then $i_{NH}^*$ is an
injection. This follows from the injection $H_T^*(Y^g)\rightarrow
H_T^*(Y^T)$ for each $g\in T$ and from
Theorem~\ref{th:productsthesame}.

\begin{corollary}\label{co:productsthesame}
  If $Y$ be a robustly equivariantly injective $T$-space, then
  the $\smile$ product on $\rous(Y)$, and the bigrading,
  can be inferred (using $i_{NH}^*$)
  from the $\star$ product and the bigrading on $\rousT$.
\end{corollary}

One way of reading the above theorem is that the cases that occur in
computing the $\smile$ product --- which $I_\lambda$ contribute to the
obstruction bundle, vs. which $I_\lambda$ are in the normal bundle hence
contribute to the $\barE_3$ pushforward --- ``cancel'' one another to
some extent when taken together, making the $\star$ product simpler
than either one considered individually.

In particular, the easy proofs of associativity and gradedness for
$\star$ imply the same for $\smile$, in the robustly equivariantly
injective case. In fact these properties hold regardless:

\begin{theorem}
  The ring $(\rous(Y),\smile)$ is bigraded and associative for any stably almost complex complex $T$-manifold $Y$.
\end{theorem}

\begin{proof}
  In the case that $Y$ is robustly equivariantly injective, we use the
  Corollary above. More generally, the proof can
  be accomplished by a case-by-case analysis
  (parallel to that in Theorem \ref{th:productsthesame}) of the bundles
  $E|_{Y^{g_1,g_2}}$ and the normal bundle $Y^{g_1,g_2}$ in
  $Y^{g_1g_2}$ over each pair $(g_1,g_2)$ and each connected component
  of $Y^{g_1,g_2}$.
\end{proof}

%% -. The non-equivariantly injective definition -- Rebecca
%%        * This will most likely get cut and moved around to other sections
%%          Rebecca is responsible for

%% 4. Relation with orbifold cohomology -- Allen
%%    relnToOrb.tex
%%        * Definition of orbifold cohomology a la Chen and Ruan, AGV
%%        * Thm: If T acts locally freely on Y, PH_T^*(Y) = H^*_orb(Y/T)

\section{Relation to orbifold cohomology of torus quotients}
\label{sec:orbifoldcohomology}
Our goal in this section is to perform the following definition chase:

\begin{theorem}\label{th:locallyfree} 
  Let $T$ act on the compact stably almost complex manifold $\locfree$
  locally freely, and let $X = \locfree / T$ be the quotient orbifold.
  Then
  $$\rous(\locfree) = H_{CR}^*(X),$$ 
  where $H_{CR}^*(X)$ is as defined in \cite{CR:orbH}.
\end{theorem}

It will become clear, as we recapitulate their definition of
$H^*_{CR}$, that we have set up our definition of $\rous$ in order to
make this tautological.

Some of the difficulty in their definition arises from the technicalities
of dealing with general orbifolds, and can be sidestepped in the case
of a global quotient. 
At one point we will need to make use of a different simplification 
of their definition, found in \cite{BCS:toricvarieties}.

\begin{proof}
\newcommand\tildeX{{\widetilde X}}
\newcommand\tildeZ{{\widetilde \locfree}}
Define 
$$
\tildeZ = \big\{ (z,g)\ |\ z\in \locfree, g \in T, g\cdot z = z\big\}\ 
\subseteq \locfree \times T.
$$
For each $z\in \locfree$, the stabilizer group is closed. By the
local freeness, each stabilizer group is discrete, hence finite.
By the compactness of $\locfree$, only finitely many stabilizer groups
occur up to conjugacy -- but since $T$ is abelian, we can omit
``up to conjugacy''. Hence only a finite set of $g$ arise this way.
Therefore the cohomology $H^*_T(\tildeZ)$ is a direct sum over these
$g\in T$, and in fact this direct sum is exactly our definition of the 
inertial cohomology:
$$ H^*_T(\tildeZ) = \rous(\locfree).$$

Since $T$ is abelian, $(z,g)\in \tildeZ$ implies $(tz,g)\in \tildeZ$ 
for all $t\in T$, so we can form the quotient by this $T$-action.
Following \cite[3.1]{CR:orbH}, we call this quotient orbifold
$\tildeX  \subseteq X\times T$.
Note that when $X$ is a manifold, $\tildeX = X \times \{1\}$.

%This space $\tildeX$ is naturally an orbifold itself, and 
We define
$$
H_{CR}^*(X) := H^*(\tildeX) \qquad\hbox{as a group.}
$$
If we work with real coefficients as in \cite{CR:orbH},
then the right hand side is 
just the ordinary cohomology of the underlying topological space.
However, we will generally prefer to use the integer cohomology of 
the classifying space of the orbifold, as in e.g. \cite{Henriques}; 
in the case at hand it means 
$$ H^*(\tildeX) = H^*_T(\tildeZ) = \rous(\locfree). $$
We're not done, though, as we still have to consider the ring structure.
(And the grading, but we leave that to the reader.)

To define the ring structure, we first need
$$ \tildeX_3 = \{ (x,g_1,g_2,g_3)\ |\ g_i\in T_x,\ g_1 g_2 g_3=1 \}, 
$$
called the {\dfn 3-multi-sector} in \cite[4.1]{CR:orbH}.
There are three natural maps
$$ e_i : \tildeX_3\to \tildeX
$$
defined by $e_i(x,g_1,g_2,g_3) = (x,g_i)$, for $i=1,2,3$. 
For each map $e_i$, we define $\overline{e_i}: \tildeX_3\to \tildeX$ by
$\overline{e_i}(x,g_1,g_2,g_3) = (x,g_i^{-1})$. 

The definition of the obstruction bundle in \cite{CR:orbH} is very complicated,
but is simplified a great deal in \cite{BCS:toricvarieties}. Let $\tildeZ_3$ be the 3-multi-sector $\tildeZ_3=\{(z,g_1,g_2,g_3)\ |\ z\in Z^{g_1,g_2},\ g_1g_2g_3=1\}.$
Let $F$ be a connected component of $\tildeZ_3$,
%mapping to the group element $g\in T$ under the projection $\tildeZ \to T$,
and $F'$ its projection to $\locfree$, considered as a component of
the space $\tildeY$ from Section \ref{sec:smileproduct}. While F is does not have an a.c.\ structure, the construction of the {\dfn obstruction bundle} $E|_F$ over $F$ of \cite[4.1]{CR:orbH} works here as well, since the normal bundle to any component of $Z^{g_1,g_2}$ in $Z$ is almost complex, and the tangent directions to $Z^{g_1,g_2}$ do not contribute to the obstruction bundle. In \cite[Proposition 6.3]{BCS:toricvarieties} the authors prove that this obstruction bundle is the quotient by T of the vector bundle F' from our Definition 3.1.%
%Then the {\dfn obstruction bundle} $E|_F$ over $F$, as defined in 
%\cite[4.1]{CR:orbH}, is proven in \cite[Proposition 6.3]{BCS:toricvarieties}
%to be the quotient by $T$ of the vector bundle over $F'$
%from our Definition \ref{de:obbundle}.
\footnote{%
  While the setting of \cite{BCS:toricvarieties} is the toric case, 
  the calculation in their Proposition 6.3 works in general.
  Obstruction theory enters \cite{CR:orbH} as the
  $H^1$ of a vector bundle constructed from the normal bundle to $F$.
  The subbundle from Definition \ref{de:obbundle} was selected out by
  asking that the sum
  $m = \Year_\lambda(g_1)+\Year_\lambda(g_2)+\Year_\lambda(g_3)$, a priori 
  either $0,1,2$, actually be $2$.
  The link provided in \cite[Proposition 6.3]{BCS:toricvarieties}
  between the two of these is to compute $H^1({\mathcal O}(-m))$ over
  ${\mathbb CP}^1$, which vanishes unless the sum is $2$.}
In \cite{CR:orbH} they consider the Euler class of this orbibundle,
as an element of $H^*(F'/T)$. In the case of a global quotient orbifold
$F'/T$, such an Euler class can instead be computed as
the equivariant Euler class of the vector bundle, 
living in the isomorphic group $H^*_T(F')$. This is exactly what we used
in the definition of $\smile$.

Let $\epsilon$ denote the sum of these Euler classes over all components,
either in $H^*_{CR}(X)$ or $\rous(\locfree)$.
Then both definitions give the product of $\alpha$ and $\beta$ as
$\overline{e_3}_*\big(e_1^*(\alpha)\cdot e_2^*(\beta) \cdot \epsilon\big)$.
\end{proof}

%% 5. Functoriality properties of preorbifold cohomology -- Allen
%%    functoriality.tex
%%        * Equivariant orbifold cohomology
%%        * If S<T, there is a restriction map PH_T --> PH_S
%%        * The orbifold cohomology is the \1 piece of PH_T^{*,*}

\section{Functoriality of inertial cohomology}
\label{sec:functoriality}
Inertial cohomology is very far from being an equivariant cohomology theory,
for much the same reasons that Chen-Ruan cohomology and quantum cohomology
fail to be properly functorial as cohomology theories.
The inertial cohomology {\em groups} are functorial:
any $T$-equivariant map $f:X \to Y$ restricts to a $T$-equivariant map
$f^t: X^t \to Y^t$ on the fixed sets by $t\in T$, and hence a map
$(f^t)^*: NH_T^{*,t}(Y)\to NH_T^{*,t}(X)$ backwards on each summand.

Since the {\em rings} depend on the stably almost complex structures --
and more specifically, the honest complex structures on normal bundles
to fixed point sets -- we will use conditions on
these to guarantee that this map $f^*$ is a ring homomorphism.
While our conditions are extremely restrictive, we have two natural
instances in which they are satisfied,
one treated here in Corollary \ref{cor:nbhds}
and one in Theorem~\ref{thm:surjectivity}.
%Section \ref{sec:surjectivity}.

\begin{proposition}\label{prop:transversehom}
  Let $i:X\into Y$ be a $T$-invariant inclusion, such that $Y$ is
  stably almost complex, and the normal bundle to $X$ in $Y$ is
  trivialized. Then $X$ is naturally stably almost complex.

  Assume also that $X$ is transverse to any $Y^t, t\in T$.
  Then the restriction map $i^*: \rous(Y)\to\rous(X)$ is a ring homomorphism.
\end{proposition}

\begin{proof}
  The trivialization of the normal bundle gives an isomorphism between
  stabilizations of the tangent bundle of $X$ and the restriction of
  the tangent bundle of $Y$. This proves the first claim.

  The transversality guarantees that for each normal bundle
  $\nu(X^t \in X)$ resp. $\nu(X^{s,t} \in X^{st})$ is the restriction of
  the corresponding normal bundle
  $\nu(Y^t \in Y)$ resp. $\nu(Y^{s,t} \in Y^{st})$, with the same logweights,
  and a simple calculation with these logweights shows the product is the same.
\end{proof}

\begin{corollary}\label{cor:nbhds}
  Let $Y$ be a stably almost complex $T$-space, and $X$ the union of
  separated $T$-invariant tubular neighborhoods of the components of $Y^T$.
  Then the obvious isomorphism of groups
  $$ \rous(X) \iso H_T^*(Y^T)\otimes_\Z \Z[T] $$
  corresponds the $\smile$ product to the $\star$ product.
  The ring homomorphism $i^*: \rous(Y)\to \rous(X)$
  composed with this isomorphism is the ring homomorphism $i_{NH}^*$
  from Section \ref{sec:smileproduct}.
\end{corollary}

\begin{proof}
  This is hardly more than a restatement of the definitions of
  $\rous(X)$ and $i_{NH}^*$. To apply Proposition \ref{prop:transversehom},
  we note that the normal bundle to $X$ in $Y$ is zero-dimensional,
  hence trivialized.
\end{proof}

Of course the best case is that $Y$ is robustly equivariantly
injective, which is exactly the statement that this restriction map
is an inclusion.

One way to think about this Corollary is the following. The ordinary
restriction map in equivariant cohomology is usually thought of as
going from $H^*_T(Y)$ to $H^*_T(Y^T)$, but could equally well go to
$H^*_T(X)$, since $X$ equivariantly deformation retracts to $Y^T$.
In the setting of inertial cohomology, by contrast, $X$ is better
than $Y^T$ through being big enough to carry the geometric information
with which we define the ring structure. Alternately, we can feed this
information in by hand, which is how we defined the $\star$ product.
In fact the idea of replacing $Y^T$ by the tubular neighborhood $X$
has already shown up in the theory of noncompact Hamiltonian cobordism
\cite{GGK}.

\newcommand{\excise}[1]{}%{$\star$\textsc{#1}$\star$}

\excise{
\todo{If the action of $T$ on $\locfree$ has a global stabilizer $A$,
  then we should relate $NH_T$ to $NH_{T/A}$. This comes up in the
  weight varieties section.}

\todo{AK: I don't think I care about any of the rest of this section,
  and would rather junk it. I'll straighten out the above at that point.}

We get more satisfying results if we think about changing the group,
such as the passage to a subgroup:

\begin{proposition}
  Let $S \into T$ be a closed subgroup, and $X$ as usual.
  Let $NH_T^{*,S}(X)$ denote the subring of $\rous(X)$ using only
  those summands $NH_T^{*,t}(X)$ with $t$ in the image of $\tau$.
  Then there is a ring homomorphism $NH_T^{*,S}(X) \to NH_S^{*,\diamond}(X)$.
  If $X$ is robustly equivariantly injective, then with rational
  coefficients this map is onto.
\end{proposition}

\begin{proof}
  As usual, the map is simple to define. For each $t\in S$,
  we have a map $ H^*_T(X^t) \to H^*_S(X^t)$ (which is rationally onto in
  the robustly equivariantly injective case). Since $S$ is a subgroup of $T$,
  any year calculation for $S$ is the same as that for $T$.
  \todo{Again, a formula would probably be more convincing}
\end{proof}

There are two ways to change the group that are particularly well behaved.

\begin{proposition}\label{prop:quotientgroup}
  Let $X,T$ be as usual, and $S$ a closed subgroup of $T$.
  \begin{itemize}
  \item Assume that $S$ acts freely on $X$.
    Then $\rous(X) \iso NH_{T/S}^{*,\diamond}(X/S)$.
  \item At the other extreme, assume that $S$ acts trivially on $X$.
    Then what?
  \end{itemize}
\end{proposition}

\todo{We want this latter case, for relating the $K$ adjoint to $K$
  simply connected in the weight varieties section}

\begin{proof}
  If $S$ acts freely on $X$, then it does so on any subset $X^t$...
\end{proof}

It would be nice to unify these two. We content ourselves with an
extended remark about the case of a {\em locally} free action by $S$.

\subsection{Toward equivariant Chen-Ruan cohomology}
We defined inertial cohomology using equivariant cohomology,
and we can extract the latter from the former:
$$ H^*_T(X) = NH_T^{*,{\bf 0}}(X). $$
This actually becomes sensible when we use
part 1 of Proposition \ref{prop:quotientgroup},
$$ \rous(X) \iso NH_{T/S}^{*,\diamond}(X/S), \qquad \hbox{$S$ acts freely;}$$
put together they say that
$$ H^*_{T/S}(X/S) = NH^{*,S/S}_{T/S}(X/S) = NH^{*,S}_T(X). $$
At first sight, this would seem to indicate that
while $H^*_T(X)$ only ``feels'' the $t={\bf 0}$ summand of $\rous(X)$,
the ring $H^*_{T/S}(X/S)$ feels all the $t\in S$ summands.
However, this is meaningless for $S$ acting freely, since these
other summands are all zero.

It seems hard to locate definitions in the literature of equivariant
cohomology of orbifolds, much less equivariant Chen-Ruan cohomology.
}

%% 6. Kirwan surjectivity -- Tara
%%    kirwan.tex
%%        * Ordinary surjectivity
%%        * Equivariant surjectivity
%%        * This will be over \Q, with comments about \Z

\section{Surjectivity for symplectic torus quotients}
\label{sec:surjectivity}
In this section we relate the inertial cohomology of
a Hamiltonian $T$-space $Y$ to the Chen-Ruan cohomology of the
symplectic reduction.  Recall that the equivariant cohomology of $Y$
surjects onto the ordinary cohomology of the reduced space $Y/\! /T$.
Our first goal is the analogue for inertial cohomology,
Theorem \ref{thm:surjectivity},
showing that inertial cohomology
surjects as a ring onto the Chen-Ruan cohomology of the reduced
space. Our second goal is to compute the kernel of this map. With it one
can express $H^*_{CR}(Y/\!/T)$ by computing $\rous(Y)$ and
quotienting by the kernel of a natural map (which we describe below).
Indeed, there is a finitely generated subring $\finitepreorb(Y)$ of $\rous(Y)$
which is sufficiently large to surject onto $H^*_{CR}(Y/\!/T)$.

Suppose that $Y$ is a Hamiltonian $T$-space, with moment map $\Phi :
Y\to \algt^*$.  If $0$ is a regular value of $\Phi$, then $T$ acts
locally freely on the level set $\Phi^{-1}(0)$, and we define the {\dfn
symplectic reduction} at $0$ to be  $Y/\! / T:= \Phi^{-1}(0)/T$.
Marsden and Weinstein showed that this is a {\em symplectic} orbifold 
(or rather, this is an easy generalization of \cite{MW:reduction}).
Kirwan used a variant of Morse theory to relate the equivariant
topology of $Y$ to the topology of $Y/\! /T$.

\begin{theorem}[\cite{K:quotients}]\label{thm:kirwan}
  Let $Y$ be a proper Hamiltonian $T$-space, with moment map
  $\Phi : Y\to \algt^*$.  Suppose that $0$ is a regular value of
  $\Phi$, and that $M^T$ has only finitely many connected components.
  Then the inclusion $\Phi^{-1}(0)\hookrightarrow Y$ induces
  \begin{equation}\label{eqn:surj}
    \kappa\ :\ H_T^*(Y;\Q) \onto H_T^*(\level;\Q) \iso H^*(Y/\! /T;\Q)
  \end{equation}
  a surjection in equivariant cohomology.  The map $\kappa$ is
  called the {\dfn Kirwan map}.
\end{theorem}

\begin{remark}
  The fact that $0$ is a regular value of $\Phi$ implies that $T$ acts
  locally freely on $\level$. This implies the isomorphism on the right
  hand side of \eqref{eqn:surj}.
\end{remark}

\begin{remark}\label{rem:surjinteger}
This theorem does hold for other coefficient rings under additional
hypotheses.  A sufficient condition is that the Atiyah-Bott
principle apply with the given coefficients, for each critical set of
$||\Phi||^2$.  For circle actions, this is merely an additional
hypothesis on the topology of the fixed point components.  For the 
action of a torus $T$, this becomes a hypothesis on the
$K$-equivariant cohomology of a critical set $C$, where
$K\subseteq T$ is the largest subtorus acting locally freely on $C$.
To apply the Atiyah-Bott principle over $\Z$, for example, it is
sufficient to assume that $H_K^*(C;\Z)$ is torsion-free. For further
details, see \cite{AtiyahBott} and \cite{TW:symplecticquotients}. 
\end{remark}

Now we turn our attention to the relationship between the inertial
cohomology of $Y$ and the Chen-Ruan cohomology of the reduced space.
As we have assumed that $0$ is a regular value of $\Phi$, the action
of $T$ on the level set $\level$ is locally free.
In particular,
this implies that the symplectic reduction is naturally an orbifold.

\begin{theorem}\label{thm:surjectivity}
  Let $Y$ be a proper Hamiltonian $T$-space, with moment map
  $\Phi : Y\to \algt^*$.  Suppose that $0$ is a regular value of $\Phi$.
  Then the inclusion $\level \hookrightarrow Y$ induces a ring homomorphism
  \begin{equation}\label{eq:orbifoldsurjectivity}
    \kappa_{NH}\ :\ \rous(Y) \to \rous(\level)
  \end{equation}
  and the latter ring is isomorphic to $H^*_{CR}(Y/\! /T)$.

Moreover, under the  assumption that $Y^T$ has only finitely many
connected components, $\kappa_{NH}$ is surjective over the rationals. 
\end{theorem}

\begin{proof}
  Since $Y$ is symplectic, its tangent bundle has a canonical $T$-invariant
  almost complex structure, up to isotopy (no ``stably'' required).
  To show that $\level$ is stably almost complex, and that $\kappa_{NH}$ is
  a ring homomorphism, we will apply Proposition \ref{prop:transversehom},
  so we establish now its two requirements. Both use the exact sequence
  $$ 0 \to T\level \into TY \onto \algt^* \to 0 $$
  which in turn depends on $0$ being a regular value.

  First, we can use the exact sequence to trivialize the normal bundle to
  $\level$, canonically up to isotopy.
  For the second, let $g\in T$, and $y \in \level \cap Y^g$.
  The component $F \subseteq Y^g$ containing $y$ is a Hamiltonian $T$-manifold
  with moment map $\Phi \circ i_F$, where $i_F$ is the inclusion $F\to Y$.
  Since $y\in \level$, its $\algt$-stabilizer is trivial, and therefore
  the differential $T(\Phi\circ i_F) : T_y F \to \algt^*$ is onto.
  By the exact sequence above, this ontoness tells us that
  $T_y F$ is transverse to $T_y \level$ inside $T_y Y$.
  Now apply Proposition \ref{prop:transversehom}.

  Since $\Phi$ was assumed proper, $\level$ is compact.
  The isomorphism of $\rous(\level)$ with $H^*_{CR}(Y/\! /T)$
  then follows immediately from Theorem~\ref{th:locallyfree}.

  For each $g\in T$, we have $(Y^g)^T = Y^T$,
  so $(Y^g)^T$ has only finitely many connected components.
  Hence we can apply ordinary Kirwan surjectivity, Theorem \ref{thm:kirwan},
  to each map $$H^*_T(Y^g; \Q) \to H^*_T(\level^g; \Q).$$
  Summing these together, we find that the map $\kappa_{NH}$ is
  surjective over the rationals.
\end{proof}

Kirwan's result gives an implicit description of the kernel of $\kappa$.
%This description is useful for understanding the kernel for
%reductions of affine space.  In that case, the kernel is the ideal
%generated by the equivariant Euler classes of the maximal co\"ordinate
%subspaces that are disjoint from $\Phi^{-1}(0)$.
%
%\todo{Is this right?  Is it
%clear?}
%
Tolman and Weitsman give an explicit description of the
kernel \cite{TW:symplecticquotients}, which will be useful
to compute the kernel for reductions of coadjoint orbits.

\begin{theorem}[Tolman-Weitsman]\label{th:TWsurjectivity}
  Let $Y$ be a compact Hamiltonian $T$-space. Let $(Y^T)_{cc}$ be the
  set of connected components of the fixed point set $Y^T$. Choose any
  $\xi\in \mathfrak{t}$ and let
  $$
  K_\xi = \{\alpha\in H_T^*(Y): \alpha|_F = 0\mbox{ for all } F\in
  (Y^T)_{cc}\mbox{ such that }\langle\Phi(F),\xi\rangle \geq 0\}.
  $$
  The kernel of the Kirwan map $\kappa$ in Equation
  (\ref{eqn:surj}) is given by the ideal
  $$\ker\kappa = \langle \cup_{\xi\in\mathfrak{t}} K_\xi\rangle.$$
\end{theorem}

%The Atiyah-Bott lemma and
The methods introduced in
\cite{TW:symplecticquotients} allow us to generalize
to the case that $Y$ is not compact, but is a proper Hamiltonian
$T$-space (its moment map has a component that is bounded from
below). This applies, for example, in the case that $Y = \C^n$ with
a proper moment map. We rephrase the theorem in this light.

We begin by asserting the existence of certain natural cohomology
classes. Let $\Phi^\xi = \langle \Phi, \xi\rangle$ be a component of
the moment map.
%, where $\xi$ is rational.
%Let $S^1\subset T$ be generated by $\xi$, and let $Y^{S^1}$ denote its fixed point set. Let $p\in Y^T$, and let $N$ be a connected component $Y^{S^1}$ containing $p$,
%(which implies $\Phi^\xi(N)$ is constant).
We follow \cite{GHJ:distinguish} and let the {\em extended stable set}
of $F$ be the set of $x\in Y$ such that there is a sequence of
critical sets $C_1,\dots,C_m\supseteq F$ such that $x$ converges to an
element of $C_1$ under the flow of $-\grad f$ (defined using a
compatible Riemannian metric), and there exist points in the negative
normal bundle of $C_j$ that converge to $C_{j+1}$ under $-\grad f$. We
say that $\alpha_F$ is a {\em Morse-Thom class} associated to $F$ if
$\alpha_F$ is homogeneous and satisfies:
\begin{enumerate}
\item $\alpha_F|_{F'} =0$ for any fixed point component $F'\in
  (Y^T)_{cc}$ not contained in the extended stable set (with respect to
  $\Phi^\xi$) of $F$ (and in particular, $\alpha_F|_{F'}=0$ if
  $\Phi^\xi(F')<\Phi^\xi(F)$).
\item $\alpha_F|_F = e_T(\nu^-F)$, the equivariant euler class of the
  negative normal bundle (with respect to $\Phi^\xi$) to $F$.
\end{enumerate}

Morse-Thom classes are not necessarily unique; there may be more
than one associated to a particular $F$.
%Fixed points are always
%contained in $Y^{S^1}$, or the critical set of $\Phi^\xi$, for any $\xi$. For sufficiently generic $\xi$,
%all the critical sets will be fixed points.

\begin{theorem}\label{th:noncompactTW}
  Let $Y$ be a proper Hamiltonian $T$-space.
%Let $\Phi^\xi$ be a component of the moment map which is bounded below.
  Let $\kappa: H_T^*(Y) \rightarrow H^*(Y/\!/T)$ be the Kirwan map of
  Equation~(\ref{eqn:surj}), where reduction is done at the regular
  value 0. The kernel of $\kappa$ is generated by
  $$
  \bigcup_{\xi\in\algt}\ \{\alpha_F\in H_T^*(Y)\mbox{ such that }
  F\in (Y^T)_{cc} \mbox{ and }\Phi^\xi(F)>0\}
  $$
  where $\alpha_F$ is any Morse-Thom class associated to
  $F$.%, and $\algt^\Q$ is the set of rational vectors in $\algt$.
\end{theorem}

We can use Theorems~\ref{th:TWsurjectivity} and
\ref{th:noncompactTW} to find the kernel of the surjection
(\ref{eq:orbifoldsurjectivity}).
\begin{corollary}\label{co:kernelpreorbtoquotient}
  Let $Y$ be a proper Hamiltonian $T$ space,
  and $\kappa_{NH}: \rous(Y)\rightarrow H_{CR}(Y/\!/T)$ the natural
  surjection. Then
  $$
  \ker\kappa_{NH}=\bigoplus_{g\in T} \ker\kappa_g,
  $$ where
  $\ker\kappa_g$ is generated those $\alpha_F\in H_T^*(Y^g)$ described
  in Theorem~\ref{th:noncompactTW}.
\end{corollary}

One immediate observation from this corollary is that, for most $g\in T$,
the entire piece $\gpiece(Y)=H_T^*(Y^g)$ is in the
kernel. This follows from the fact that, for generic $g$, $Y^g = Y^T$
(and misses $\level$). Indeed, the only values of $g\in T$ such that
$\gpiece(Y)$ is not contained in $\ker\kappa_{NH}$ are those such that
$Y^g$ has an effective $T$ action. In other words, they are finite
stabilizers. We find a smaller ring which surjects onto
$H^*_{CR}(Y/\!/T)$, by excluding those $\gpiece(Y)$ in $\rous(Y)$ such
that $g$ is not a finite stabilizer.

%\todo{This is the cut-and-pasted finite stabilizer stuff.  I need to
%rewrite it.}

%We note that because $\kappa_{NH}$ is a ring map, we may compute the
%Chen-Ruan cohomology {\em ring} by computing the inertial
%cohomology ring, and determining the kernel of the map
%$\kappa_{NH}$. In general, the inertial cohomology ring is rather
%unwieldy; however, most of the inertial cohomology ring is in the
%kernel of $\kappa_{NH}$.  Of particular importance is the finite
%subring of $\rous(Y)$ obtained by restricting our attention to those
%pieces $\gpiece(Y)$ for which $g$ occurs in the group generated by
%the finite stabilizers. We introduce some notation:

\begin{definition}
  An element $g\in T$ is a {\dfn finite stabilizer} if there exists a
  point $y\in Y$ with $\stab(y)$ finite and $g\in \stab(y)$.  We let
  $\Gamma$ denote the group generated by all finite stabilizers of $Y$
  in $T$.
  %Let $\Gamma_y=\stab(y)$ when the stabilizer is finite, and
  %otherwise let $\Gamma_y=\emptyset$. What was this thing anyway?
  We assume that $\Gamma$ is finite, as is automatic if $Y$ is of
  finite type.
\end{definition}

\begin{remark}\label{re:rationalage}
  If $g\in \Gamma$, then $\Year_{\lambda_j}(g), j = 1,\dots, n$ are
  rational numbers. In other words, the grading restricted to
  $\finitepreorb(Y)$ is rational. This accounts for the rational
  grading on the Chen-Ruan cohomology of the quotient space. See
  \cite{FG:globalquotients} and \cite{R:age}.
\end{remark}

\begin{lemma}\label{lem:finstabs}
  An element $g\in T$ is a finite stabilizer on $Y$ if and only if
 there exists $p \in Y^T$ such that the weights $\lambda \in \hatT$ of $T_p Y$
  with logweight $\Year_\lambda(g) = 0$ linearly span the
  weight lattice $\hatT$ (over $\Q$).
\end{lemma}

\begin{proof}
  First, we note that an element $g$ is a finite stabilizer if and only if
  $Y^g$ contains a component on which the generic $T$-stabilizer is finite.
  Equivalently, the $\mathfrak t$-stabilizer should be trivial.

  Let $F$ be a component of $Y^g$ on which the stabilizer has minimum
  dimension, and $p \in F^T$. Then the generic stabilizer on $F$ is the
  same as the generic stabilizer on $T_p F$, which is the intersection
  of the kernels of the weights $\lambda$ on $T_p F$. For this intersection
  to be zero, then dually, the weights $\lambda$ should span $\hatT$.
\end{proof}

\begin{definition}
  The {\dfn $\finitepiece$} $\finitepreorb(Y)$ of $\rous(Y)$ is given
  as an $H^*_T(pt)$-module by
  \begin{align*}
    \finitepreorb(Y) :&= \bigoplus_{g\in \Gamma}\gpiece(Y).
  \end{align*}
  It follows from Remark~\ref{re:whereproductlives} that
$(\finitepreorb(Y),\smile)$ is s subring of $\rous(Y)$.
\end{definition}

%We note that
%for any $g\not\in \Gamma$, the vector space $\gpiece(Y)$ is in the
%kernel of the surjection (\ref{eq:surjection}).  This is because there
%is no reduction that has the element $g$ as an element of one of the
%local stabilizing groups.

\noindent The following corollary immediately follows.
\begin{corollary}\label{co:finiteorbifoldsurjectivity}
  Let $Y$ be a proper Hamiltonian $T$-space. Suppose that $Y^T$ has
  only finitely many connected components. By abuse of notation, we
  write $\kappa_{NH}$ for $\kappa_{NH}|_{\finitepreorb(Y)}$.  Then
\begin{equation}\label{eq:finiteorbifoldsurjectivity}
\xymatrix{
\kappa_{NH}\ :\ \finitepreorb(Y;\Q) \ar@{->>}[r] &
\rous(\level;\Q) \iso H^*_{CR}(Y/\! /T;\Q) },
\end{equation}
where the reduction is taken at any regular value
of the moment map.  As before,
$$
\ker\kappa_{NH}=\bigoplus_{g\in T} \ker\kappa_g,
$$ where
$\ker\kappa_g$ is generated those $\alpha_F\in H_T^*(Y^g)$ described
in Theorem~\ref{th:noncompactTW}.
\end{corollary}
Thus $H_{CR}^*(Y/\! /T)$ may be computed by finding the (finitely
generated) ring $\finitepreorb(Y)$ and quotienting by the (finitely
generated) $\ker\kappa_{NH}$.  We show several such
computations in Sections~\ref{sec:singularities} and
\ref{sec:toricvarieties}.

%% 7. EGs: Hamiltonian spaces with isolated fixed points -- Rebecca
%%    EGham.tex

\section{A graphical view for the product on Hamiltonian T-spaces \\
  and the inertial surjection (\ref{eq:surjection})}
\label{sec:multiplication}
In this section we assume that $Y$ has isolated fixed point set
$Y^T$. The most commonly studied examples are toric varieties and
flag manifolds, but we will even find interest in the case of $Y$ a
vector space with $T$ acting linearly. Very shortly we will also
require $Y$ to be a proper Hamiltonian $T$-space (in particular,
robustly equivariantly injective).

%We give yet another description of the product that allows one to multiply
%classes by restricting to fixed points.
% However, rather than
% calculating the experience of line bundles under several group elements,
%This method uses a new restriction map $res$ that assigns
%to each class, for every fixed point component $F$ and every $g\in T$,
%a product of Euler classes raised to {\em real} powers.

%\subsection{Multiplication on the fixed point set}

Using the standard ring structure on $\RoUST$
(i.e. {\em not} our multiplication $\star$),
the natural restriction $\RoUS(Y)\longrightarrow \RoUST$
obtained by restricting to the fixed point set on each piece
is not usually a ring homomorphism. To make it one, we had to invent the
$\star$ product on $\RoUST$, which twisted the multiplication using logweights.

In this section we will work the logweights into the {\em homomorphism},
rather than the multiplication on the target, giving yet another
description of the multiplication. To do so, though,
we will have to enlarge our base ring.

\renewcommand\H{{\mathcal H}}
\subsection{The base ring $\H_T$ and a new restriction map}
Recall that $H^*_T(pt)$ is naturally isomorphic (over $\Z$)
to the symmetric algebra on its degree $2$ part,
the weight lattice $\hatT$ of the torus $T$.
Define the commutative $H^*_T(pt)$-algebra
$\H_T$ by
$$ \H_T := \Z\big[ \{w^r : w\in \hatT, r\in \R_{>0}\} \big] \bigg/
        \big\langle w^1 + v^1 = (w+v)^1,\ (w^r) (w^s) = w^{r+s} \big\rangle, $$
in which we have included all positive real powers of our generators
$H^2_T(pt)$. There is an evident inclusion of $H^*_T(pt)$ into $\H_T$
induced from $w\mapsto w^1$, and a grading on $\H_T$, where $\deg w^r = 2r$.
It seems worthy of note that including real powers into this ring has
not rationalized it; in particular $\H_T^0$ is still just $\Z$.

For $\alpha \in NH^{*,g}(Y)$ (i.e. of pure degree in the second component),
and $p \in Y^T$,
we define the `restriction' map
\begin{equation}\label{eq:resatp}
 res(\alpha)|_p := \alpha|_p
        \prod_{\lambda \in \hatT}
        \lambda^{\dim I_\lambda\ \Year_\lambda(g)} \qquad\in \H_T
        \end{equation}
where $I_\lambda$ is the $\lambda$ weight space of $T_p Y$.
Summing these maps together, each tensored with $g \in \Z[T]$ to record
the $T$-grading, we get a map
$$ res: \RoUS(Y) \to \bigg( \bigoplus_{p\in Y^T} \H_T \bigg) \tensor \Z[T] $$
which will take the place of $i_{NH}^*$ from Section \ref{sec:starproduct}.

\begin{theorem}\label{th:restrictionproduct}
  Let $a,b\in \RoUS(Y)$, and let $Y$ have isolated fixed points.
  Then
  $res: \RoUS(Y) \to \bigg(\bigoplus_{p\in Y^T} \H_T \bigg) \tensor \Z[T]$
  is a graded ring homomorphism, taking $\smile$ to the ordinary product.
\end{theorem}

\begin{proof}
  It is enough to check for $a\in NH_T^{*,g_1}(Y), b\in NH_T^{*,g_2}(Y)$.
  Let $p\in Y^T$. Then
  \begin{align*}
    res(a\smile b)|_p
    &=  (a\smile b)|_p
    \prod_{\lambda} \lambda^{\dim I_\lambda\ \Year_\lambda(g_1 g_2)} \\
    &= a|_p b|_p
    \prod_{\lambda}
    \lambda^{\dim I_\lambda \big( \Year^p_\lambda(g_1)+\Year^p_\lambda(g_2)
      -\Year^p_\lambda(g_1 g_2) \big)} %\in  NH_T^{*,g_1\cdots g_n}(p)
    \prod_{\lambda} \lambda^{\dim I_\lambda\ \Year_\lambda(g_1 g_2)} \\
    &= a|_p b|_p
    \prod_{\lambda}
    \lambda^{\dim I_\lambda \big( \Year^p_\lambda(g_1)+\Year^p_\lambda(g_2)
       \big)} \\
    &= \big( a|_p
    \prod_{\lambda}     \lambda^{\dim I_\lambda \Year^p_\lambda(g_1)} \big)
    \big( b|_p
    \prod_{\lambda}    \lambda^{\dim I_\lambda \Year^p_\lambda(g_2)} \big) \\
    &= (res\ a)|_p (res\ b)|_p.
  \end{align*}
  To get from the first to the second line, we used Theorem
  \ref{th:productsthesame}.

  To check the grading, we need to assume $a$ is of pure degree in
  $\RoUS(Y)$, for example if
  $a \in NH_T^{*,g}(Y) = \oplus_{F \subseteq Y^g} H^*_T(F)$ is actually
  a homogeneous element of $H^*_T(F)$ for some component $F$ of $Y^g$.
  Then for $p$ any element of $F^T$,
  $$ \deg res\ a
  = \deg_{\H_T} (res\ a)|_p
  = \deg_{\H_T} \big( \alpha|_p
  \prod_{\lambda \in \hatT} \lambda^{\dim I_\lambda\ \Year_\lambda(g)} \big)
  = (\deg \alpha_p) + \sum_\lambda 2 \dim I_\lambda\ \Year_\lambda(g) $$
  which is exactly the age-shifted definition we gave for the grading on
  the $H^*_T(F)$ component of $NH^{*,g}_T(Y)$.
\end{proof}

\begin{remark}
  The kernel of $res$ is the same as the kernel of $i_{NH}^*$.
  In particular, if $Y$ is robustly equivariantly injective,
  we can use Theorem \ref{th:restrictionproduct} to compute
  the $\smile$ product.
\end{remark}

\subsection{A pictorial description of the product}

In this subsection we assume that $Y$ is a proper Hamiltonian $T$-space,
as this will allow us to read off the finite stabilizers from information
that is often recorded with the {\dfn moment polyhedron} $\Phi(Y)$.
We recall the basic facts we will need from the geometry of moment maps,
as can be found in e.g. \cite{GS:convexity}. The pictorial description
extends that used in ordinary equivariant cohomology, as detailed in
e.g.  \cite{HHH}.
%%%\cite[Appendix]{HHH}.

The image $\Phi(Y) \subseteq \mathfrak{t}^*$ of $Y$ under the moment
map $\Phi$ is a convex polyhedron (possibly unbounded), and when $Y$
is compact, it is the convex hull of the finite set $\Phi(Y^T)$,
then called the {\dfn moment polytope}. For $p\in Y^T$, and
$\lambda$ a weight of $T_p Y$, the component $F$ of $Y^{\ker
\lambda}$ containing the point $p$ is itself a proper Hamiltonian
$T$-space. Its moment polyhedron $\Phi(F)$ is an {\em interval}
inside $\Phi(Y)$ with one end at $\Phi(p)$, continuing in the
direction $\lambda$. (To think of $\lambda \in \hatT$ as a vector in
$\mathfrak{t}^*$, we are using the natural embedding $\hatT \to
\mathfrak{t}^*$.) When we draw moment polyhedra, we will always
superimpose these intervals upon them, which include the edges of
the polyhedron $\Phi(Y)$. Therefore, from the picture alone, one can
almost determine the weights of the $T$-action on $T_p Y$ -- but only up
to positive scaling and multiplicity. {\em We will assume we know
the actual weights.}

\begin{figure}[h]
\begin{center}
  \epsfig{figure=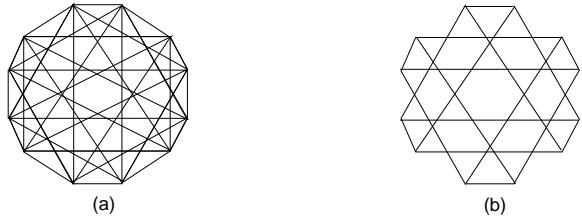,width=3.0in}
\end{center}
\caption{(a) The moment map image of a generic coadjoint orbit of $G_2$,
  and (b) the image of the fixed point set of a special order $3$ element
  of $T\subset G_2$.}
 \label{fig:G2whole}
\end{figure}

Given $g\in T$, the moment map image $\Phi(Y^g)$ may not be convex,
since $Y^g$ is not necessarily connected. The moment map image for a
generic coadjoint orbit $\mathcal{O}$ of $G_2$ under the action of
its maximal torus is shown in Figure~\ref{fig:G2whole}(a). Let $T$
be a maximal torus of $G_2$, and let $g\in T$ be an element of order
3 that fixes two copies of a generic coadjoint orbit of $SU(3)$
inside. The image of $\mathcal{O}^g$ can be seen in
Figure~\ref{fig:G2whole}(b).

%%%%%%%%%%%%%%%%%%%%%%%%%As in \cite[Appendix]{HHH},
%%%%%%%%%%%%%%%
There is a pictorial way to represent elements of $\H_T$ (and later,
$\rous(Y)$) using the fact that $\mathfrak t^*$ plays two roles;
it is the home of the moment polytope $\Phi(Y)$,
and is also the generators of $H_T^*(pt)$.
Each monomial $w_1^{r_1}w_2^{r_2}\cdots w_k^{r_k}\in \H_T$ may be
%represented pictorially
{\dfn drawn} by drawing the vectors $w_1,\dots, w_k$
and labeling each $w_i$ with the positive real number $r_i$. A sum of such monomials may be drawn as formal sum of such vector drawings, one for each monomial.

As noted previously, the {\em res} map is injective
for Hamiltonian spaces.
So to draw a class $\alpha_g \in \gpiece(Y)$,
we can consider its image under the {\em res} map.

\begin{figure}[h]
  \begin{center}
    \epsfig{figure=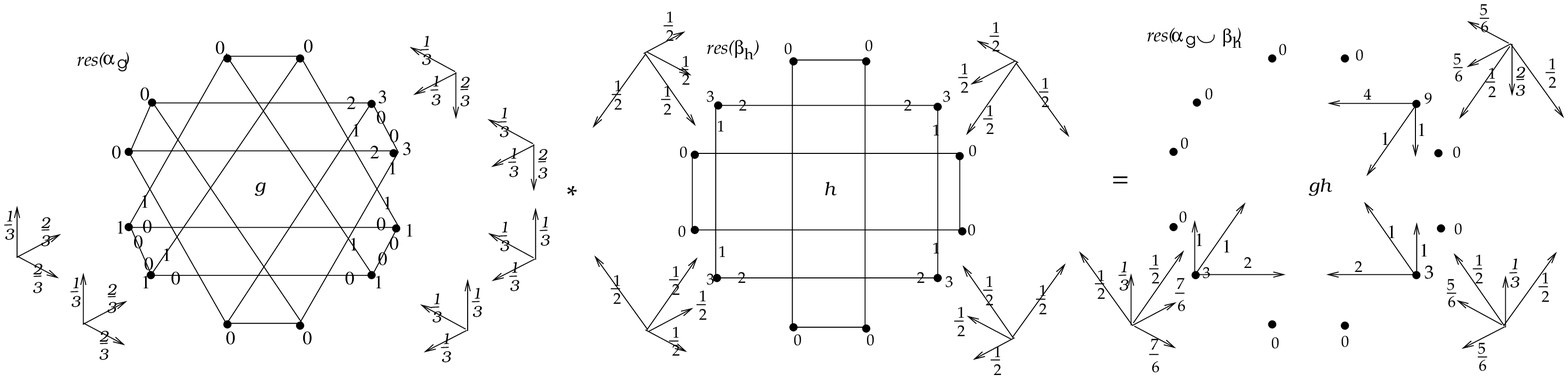,width=6.5in}
  \end{center}
  \caption{The product of two diagrams.}
  \label{fig:productofdiagrams}
\end{figure}

Let $p$ be a fixed point, which we assume to be isolated.
%By Equation~\ref{eq:resatp}, $res(\alpha_g)|_p \in \H_T$.
Note that $res(\alpha_g)|_p$ may be
expressed as a product of two elements of $\H_T$:
\begin{itemize}
\item $\prod_{\lambda \in \hatT} \lambda^{\dim I_\lambda\ \Year_\lambda(g)}$,
  which depends on $g$ but not on $\alpha$.
\item $\alpha_g|_p$.
\end{itemize}
Suppose we know $(\alpha_g)|_p$ explicitly, and it can be written as a monomial $c\cdot w_1^{r_1}\cdots w_k^{r_k}$, with $r_i\in \Z$. We draw $w_1^{r_1}\cdots w_k^{r_k}$ at $\Phi(p)$, and any coefficient $c$ as a number at $\Phi(p)$. We draw the product $\prod_{\lambda \in \hatT} \lambda^{\dim I_\lambda\ \Year_\lambda(g)}$ near the moment image $\Phi(X^g)$, close to $\Phi(p)$. This has the advantage of separating these two different pieces of the computation of $res(\alpha_g)$. If $(\alpha_g)|_p$ is a sum of such monomials, then $res(\alpha_g)|_p$ is a sum of these labeled-vector drawings, one for each monomial. For later purposes, we caption each labeled-vector drawing
with the element $g$ as well. If $\alpha_g|_p$ is not given explicitly, it may be more convenient to draw the vectors in slightly different positions; see Figure~\ref{fig:productofdiagrams}.

There is a slight annoyance if two fixed points $p$ and $q$ have
$\Phi(p) = \Phi(q)$, in that we have to move the picture of one of them.
(This is not a mathematical objection, just a practical one.)

\begin{figure}[h]
\begin{center}
\epsfig{figure=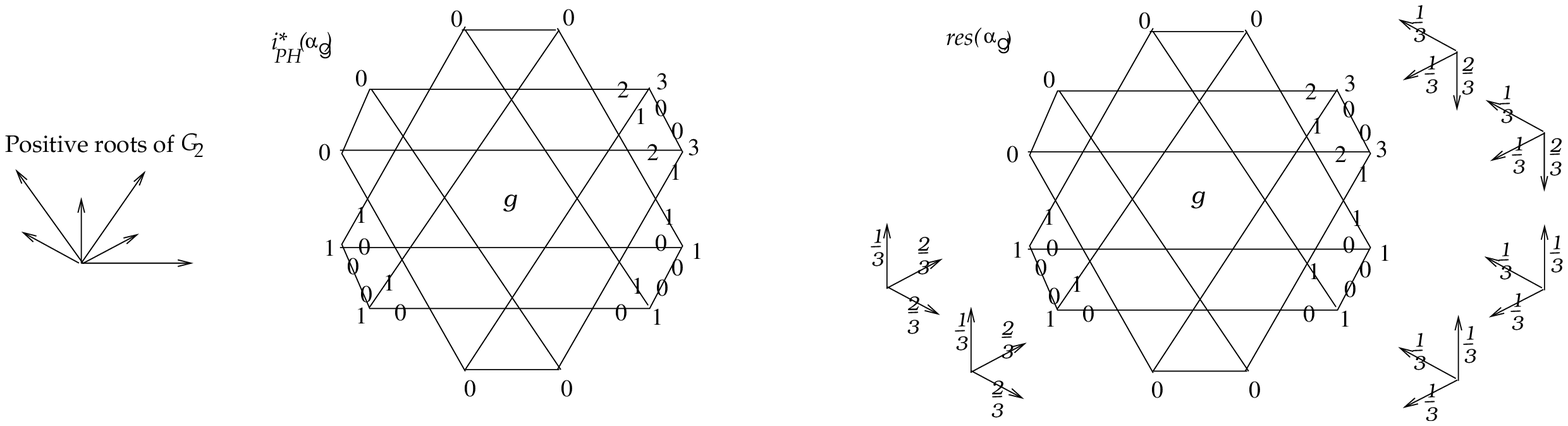,width=6.5in}
\end{center}
\caption{The restrictions $i_{NH}^*(\alpha_g)$ and  $res(\alpha_g)$,
  drawn on the picture $\mathcal{O}^g$, the fixed point set under $g$
  of a generic coadjoint orbit of $G_2$, where $g$ is an element of
  order 3 fixing two coadjoint orbits of $SU(3)$.}
\label{fig:alphagandresalphag}
\end{figure}

We may now draw $res(\alpha_g)$ by associating a sum of diagrams to
each $p\in Y^T$. In particular, if $res(\alpha_g)|_p$ is a
monomial for every $p$, then one diagram suffices to represent
$\alpha_g$. In Figure~\ref{fig:alphagandresalphag} we draw the
picture corresponding to $i_{NH}^*(\alpha_g)=\alpha_g|_{Y^T}$, and that
corresponding to $res(\alpha_g)$.

While each $\alpha_g$ may be written as a class on $H_T^*(Y^g)$,
$res(\alpha_g)$ may not be drawn on the moment map for $Y^g$ (using
only the weights of $T$ on $T_pY^g$ for each $p$).
At any fixed $p$, the vectors $\lambda$ occurring in the term
$\prod_{\lambda \in \hatT} \lambda^{\dim I_\lambda\ \Year_\lambda(g)}$
of Equation~\ref{eq:resatp} point out of $Y^g$; they are by definition
those $\lambda$ occurring in $T_pY$ whose logweights
$\Year_\lambda(g)$ are {\em not} 0.

Multiplication of two classes
$\alpha_g\in \gpiece(Y)$ and $\beta_h\in NH_T^{\ast,h}(Y)$
is easy in this pictorial calculus.
The product of classes is performed pointwise,
and involves only the product structure on $\H_T$
(with no additional factors such as those introduced by the $\star$
product, since they've been worked into $res$). The $\Z[T]$ factor
is only there to remember that $\alpha_g \beta_h$ lives in
$NH_T^{\ast,gh}(Y)$.

By distributivity, it is enough to
treat the case that each $\alpha_g|_p$ or $\beta_h|_p$ is a monomial.
The product of a diagram labeled by $g$ and one labeled by $h$ is
labeled by $gh$. The label (exponent) on a vector $\lambda$ at $p$ in the
product is the sum of the labels at $p$ in the $g$-diagram and the
$h$-diagram.

For example, let $g$ be an order 3 element in the maximal torus of
$G_2$ fixing two copies of $SU(3)/T$ in $G_2/T$, and $h$ an order 2
element fixing three copies of $SO(4)/T$. Two elements
$\alpha_g$ and $\beta_h$ and their product are described by the
diagrams in Figure~\ref{fig:productofdiagrams}.

Note that the final picture is obtained by adding the labels of each of the vectors, but multiplying the coefficients at $\Phi(p)$ for each fixed $p$. The result is a pictorial representation of $res(\alpha_g\smile\beta_h)$; however, it is not separated into $(\alpha_g\smile\beta_h)|_p$ and $\prod \lambda^{\dim I_\lambda a_\lambda(gh)}$. The monomial drawn at $\Phi(p)$ is $\alpha_g|_p\beta_h|_p$, and the vectors near $\Phi(p)$ represent $\prod \lambda^{\dim I_\lambda (a_\lambda(g)+a_\lambda(h))}$. These product of these pieces are the same by Theorem~\ref{th:restrictionproduct}.

\subsection{Finding the finite stabilizers}\label{se:finitestab}

It is clear at this point that calculating $\finitepreorb(Y)$ instead
of $\rous(Y)$ has an appeal: there are finitely many labeled diagrams
such that all elements of $\finitepreorb(Y)$ may be expressed as
$H^*_T(pt)$-linear combinations of these diagrams.
And as we noted in Corollary~\ref{co:finiteorbifoldsurjectivity}, the ring $\finitepreorb(Y)$ is large
enough to surject on to the Chen-Ruan cohomology of the reduced space
$Y/\!/T$.

To determine the finite stabilizers, we first need a way to picture an element
$g$ of $T$. By Pontrjagin duality, we see that $T \iso \Hom(\hatT,U(1))$,
so we can reconstruct $g$ from the function labeling each point
$\lambda\in \hatT$ by the logweight $\Year_\lambda(g)$. We can thus picture
$g$ as a labeling of a generating set of $\hatT$ by elements of $[0,1)$,
and require that the logweights come from a homomorphism.

Lemma \ref{lem:finstabs} said that an element $g$ is a finite stabilizer if
there exists $p\in X^T$ such that the set
$$ \{\lambda \hbox{ is a weight of $T_p X$, and annihilates $g$} \} $$
is big enough to $\Q$-span $\mathfrak t^*$.
If we assume that $T$ acts faithfully on $X$,
then the union over $p\in X^T$ of the weights at $T_p X$ will span $\hatT$.

\begin{figure}[h]
  \begin{center}
    \epsfig{figure=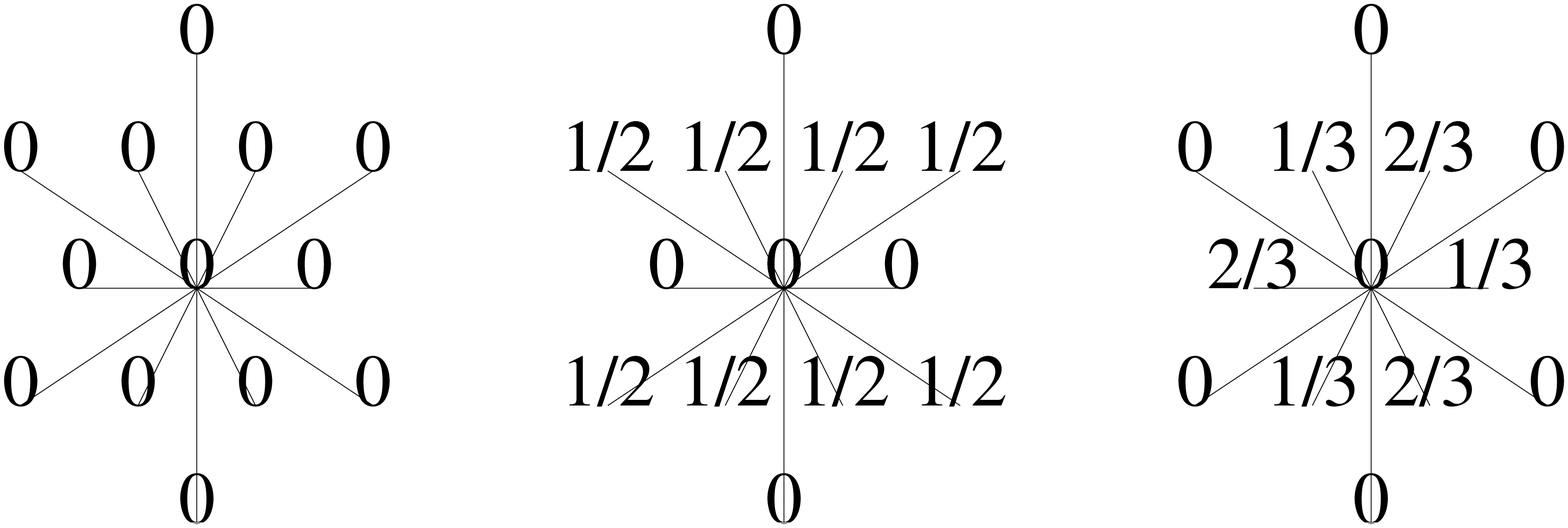,width=3in}
  \end{center}
  \caption{Three finite stabilizers of the action on a coadjoint orbit
    of $G_2$, pictured by their logweights on the root system.}
  \label{fig:G2finstab}
\end{figure}

We illustrate this technology to find the finite stabilizers in the
example of $X$ a coadjoint orbit of $G_2$, where the union of the
weights at all fixed points is exactly the root system of $G_2$.
There are three finite stabilizers, up to rotation and reflection (the
action of the Weyl group of $G_2$), pictured in Figure \ref{fig:G2finstab}.

\subsection{Finding the kernel of Equation \ref{eq:surjection} using pictures}

Before proceeding to a first completely worked example, we illustrate
the use of pictures in describing the kernel of the Kirwan map
$$ \kappa_g: H_T^*(Y^g) \longrightarrow H_T^*(Y^g/\!/T),
$$
and the analogous kernel of the inertial Kirwan map
$$ \kappa_{NH}:\rous(Y) \longrightarrow H_{CR}(Y/\!/T).
$$
Recall Corollary~\ref{co:kernelpreorbtoquotient} states that $\kappa_{NH}$
is generated by elements in the kernel of $\kappa_g$ for any $g\in T$.

By the Tolman-Weitsman theorem (Theorems~\ref{th:TWsurjectivity} and \ref{th:noncompactTW}), the
kernel of $\kappa_g$ is generated by classes $\alpha\in H_T^*(Y^g)$
that satisfy the property that
\begin{center}
  \begin{minipage}[t]{.8\linewidth}
    there exists $\xi\in \mathfrak{t}$ such that $\alpha|_p=0$ for all
    $p$ with $\langle \xi, \Phi(p)\rangle \geq 0$.
  \end{minipage}
\end{center}
Then the kernel of $\kappa_{NH}$ is generated by classes $\alpha\in
\rous(Y)$ such that $\alpha\in \gpiece(Y)$ and satisfies this
property, for some $g$.
% For reduction at a point $\mu$ different from
% 0, replace the value 0 with $\langle \xi,\mu\rangle$.
% AK: This 0 is ambiguous above, and it breaks the flow to say enough
% to make it sensical. Those few who care will know what to do.

\begin{figure}[h]
  \epsfig{figure=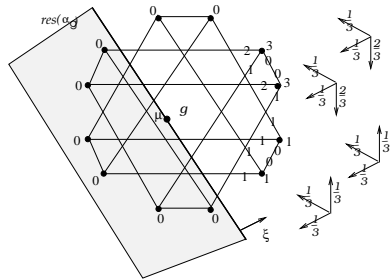,width=2in}
  \caption{A class in the kernel of map from inertial cohomology
    to Chen-Ruan cohomology of the reduction.}\label{fig:kernelelement}
\end{figure}

In Figure~\ref{fig:kernelelement} we show an example of an element in
the kernel of the map
$$
\rous(G_2/T)\longrightarrow H_{CR}(X),
$$
where $X$ is the orbifold obtained by symplectic reduction of this
coadjoint orbit by $T$ at 0.
For those familiar with equivariant cohomology of coadjoint orbits (or
more generally, of GKM spaces), one might notice that the labeling on
$\Phi(Y^g)$ must be compatible in some sense: the part drawn on
$\Phi(Y^g)$ is the class restricted to $(Y^g)^T$ in equivariant
cohomology.

To obtain the kernel of $\kappa_{NH}$, one must take every class in
$\gpiece(Y)$ that has this property for some $\xi\in \mathfrak{t}$,
and then to do this for every picture, as $g$ varies in $T$.

\subsection{A toric example: $\C_1 \oplus \C_1 \oplus \C_3$}

Let $T=S^1$, acting on $\C^3$ with weights $1$, $1$, and $3$. We
will call the three weight lines $\C_1, \C_1', \C_3$. Then the
finite stabilizers are $g_s=\exp(2\pi i s)$,
$s=0,\frac{1}{3},\frac{2}{3}$. Encoded as functions on the weight
lattice $\Z$, they are
$$ \ldots, 0, 0, {\bf 0}, 0, 0, \ldots $$
$$\ldots, {1/3}, {2/3}, {\bf 0}, {1/3}, {2/3}, \ldots$$
$$\ldots, {2/3}, {1/3}, {\bf 0}, {2/3}, {1/3}, \ldots$$
In this case, they don't only generate but are already equal to
the subgroup $\Gamma \leq T$.
The fixed point sets for these group elements $g_0, g_{1/3}, g_{2/3}$ are
$$ \C_1 \oplus \C_1 \oplus \C_3, \quad \C_3, \quad \C_3 $$
respectively, so each $\gpiece(\C^3)$ is free of rank $1$
over $H^*_T(pt) = \Z[u]$.

We calculate in detail the $res$ map applied to the generator of
$NH^{*,1/3}(\C^3)$ at the one fixed point, $\{0\}$. This is a product
over $\C_1, \C_1', \C_3$ of $u$ raised to the logweight power of
$g_{1/3}$, respectively $1/3$, $1/3$, $0$.  Then we tensor with
$g_{1/3}\in T$ to keep track of the $T$ grading.  The result is
$u^{2/3} \tensor g_{1/3}$.

In all, the three generators have $res|_{\{0\}}$ of
$$ 1\tensor g_0, \quad u^{2/3} \tensor g_{1/3},\quad u^{4/3} \tensor g_{2/3}.$$
If we call these $1,a,b$, then $a^2 = b$, $a^3 = u^2$. So
$$
NH^{*,\Gamma}(\C^3) =\Z[u,a,b]/\langle a^2-b, a^3-u^2\rangle=
\Z[u,a] / \langle a^3 - u^2 \rangle $$
where the bidegree of $a$ is
$(4/3, g_{1/3})$. (Recall that we get an extra factor of $2 = \dim_\R
\C$ through working with cohomology rather than Chow rings.)

It is left to calculate the Chen-Ruan cohomology of the symplectic
quotient at a regular value. Let us reduce at 1, and let $X$ be the
orbifold quotient, i.e.
$$
X = \{|z_1|^2 + |z_2|^2+ 3|z_3|^2 = 1\}/S^1.
$$
According to the noncompact version of the Tolman-Weitsman theorem
(Theorem~\ref{th:noncompactTW}), the kernel is generated by
classes $\alpha\in \gpiece(\C^3)$ (for each $g\in T$) whose
restriction to $\{0\}$ is a multiple of the Euler class of the
(negative) normal bundle to $\{0\}$. As discussed in Section~\ref{se:integers},
these classes will generate the kernel also when the inertial cohomology
is taken with $\Z$ coefficients. As a module over $H_T^*(pt)$, we have
$$
\finitepreorb(\C^3) = H_T^*(\C_1\oplus \C_1\oplus \C_3)\oplus
H_T^*(\C_3)\oplus H_T^*(\C_3)
$$
The equivariant Euler class of $\{0\}$ in the first piece is
$u\cdot u\cdot 3u = 3u^3$. In the second piece the equivariant Euler
class of $\{0\}$ is $3ua$, and for the third piece is $3ub = 3ua^2$.
We thus obtain
$$
H^*_{CR}(X) = \Z[u^{(2)},a^{(4/3)}] / \langle a^3 - u^2, 3u^3, 3ua \rangle
$$
where the superscripts indicate the degree.

If we drop the generator $a$, and rationalize, we get the ordinary
cohomology $H^*(X;\Q) = \Q[u]/\langle u^3\rangle$ of the coarse moduli space.

% Rationally, this is
% $$ H_{CR}^0 = \Q,\ H_{CR}^{4/3} = \Q a,\ H_{CR}^2 = \Q u,
% H_{CR}^{8/3} = \Q a^2,\ H_{CR}^4 = Q u^2. $$

%% 8. EGs: Weight varieties -- Allen
%%    EGweight.tex
%%        * What finite stabilizers occur for T acting on a coadjoint orbit
%%        * T-reductions of B_2, G_2 worked out explicitly

\section{Flag manifolds and weight varieties}
\label{sec:singularities}
In this section we study the example of $Y = K/T$, 
called a {\dfn generalized flag manifold}, where $K$ is a 
compact, connected Lie group and $T$ is a maximal torus thereof.

This example is already well handled by the techniques of the
last section, as it is Hamiltonian (which we will go over in a moment)
with isolated fixed points, and all the fixed points map 
to different places under the moment map. 
The main result of this section is then an efficient calculation of
the finite stabilizers.  One interesting corollary is that if $K$ is a
classical group (and only if), the obstruction bundles are all trivial.

%which acts on the left. 
The standard notation we need from Lie theory is the
normalizer of the torus $N(T)$, the Weyl group $W := N(T)/T$,
and the centralizer $C_K(k)$ of an element $k \in K$. 
%We will also have need of the simply-connected cover $\widetilde K$ of $K$.  
This space $K/T$ has a left action of $K$ and hence of $T$, and a right
action of $W$. It has a family of symplectic structures, one for each
generic orbit $K\cdot \lambda$ on the dual $\lie{k}^*$ of the Lie algebra
of $K$.
The moment map $\Phi$ is the projection $\lie{k}^*\to\lie{t}^*$ transpose
to the inclusion of Lie algebras. Using the Killing form, we can and will
regard the basepoint $\lambda$ as an element of $\lie{t}^* \into \lie{k}^*$.

A {\dfn weight variety} \cite{AllenThesis,Go:thesis} is the
symplectic quotient of a coadjoint orbit $K\cdot \lambda$ by the
maximal torus $T\leq K$.  These turn out to be smooth (for reductions
at regular values) for $K=SU(n)$, but are orbifolds for other $K$,
as we will explain after Proposition \ref{prop:specialclasses}.

We will need a few standard facts about such $K$:
\begin{itemize}
\item Every element of $K$ is conjugate to some element of $T$.
\item If two elements of $T$ are $K$-conjugate, they are already
  conjugate by $N(T)$.
\item The center $Z(K)$ is contained in $T$. 
\item Any two maximal tori in $K$ are conjugate.
\end{itemize}
The group $K$ is {\dfn semisimple} if its center $Z(K)$
is finite, or equivalently, if the center of its Lie algebra is trivial.

\begin{lemma}\label{lem:genericstab}
  Let $K$ be a compact connected Lie group, and $T$ a maximal torus,
  so $T$ contains the center $Z(K)$.
  Then the generic stabilizer of $T$ acting on $K/T$ is $Z(K)$.
\end{lemma}

In particular, unless $K$ is semisimple,
there are no finite stabilizers at all.

\begin{proof}
  Since $Z(K) \leq T$, for all $z\in Z(K), g\in K$ we have
  $$ z g T = g z T = g T. $$

  Conversely, let $s\in G$ stabilize every point of $K/T$, 
  so $\forall k\in K$, $s k T = k T$, hence $s\in k T k^{-1}$,
  so $s$ commutes with $k T k^{-1}$. But the union over $k\in K$ of
  the tori $k T k^{-1}$ is all of $K$, since every element can be
  conjugated into $T$. Hence $s$ commutes with all of $K$.
\end{proof}

\subsection{The finite stabilizers in $T$ on $K/T$}

We are now ready to determine, following \cite{AllenThesis}, which
$T$-stabilizers occur on $K/T$, and their fixed points.

\begin{lemma}\label{lem:stabs}
  \begin{itemize}
  \item
    Let $k\in K$, and $kT$ the corresponding point in $K/T$.
    Then $kT$ is stabilized by $t\in T$ if and only if $k \in C_K(t) N(T)$.
  \item 
    Let $C_K(t)^0$ denote the identity component of $C_K(t)$, and $W_t$
    denote the Weyl group of $C_K(t)^0$ (with respect to the same 
    maximal torus, $T$). 
    Then $C_K(t) N(T) = C_K(t)^0 N(T)$.
    Each component of $C_K(t) N(T) / T$ 
    is isomorphic to the smaller flag manifold $C_K(t)^0/T$,
    and the components are indexed by the cosets $W_t \backslash W$.
  \item 
    An element $t\in T$ occurs as a finite stabilizer if and only
    if the identity component $C_K(t)^0$ of $C_K(t)$ is semisimple,
    and necessarily of the same rank as $K$.
  \end{itemize}
\end{lemma}

\newcommand\Iff{\quad\Longleftrightarrow\quad}

\begin{proof}
  To start off the first claim, 
  $$ t kT = kT \Iff k^{-1}tkT = T \Iff k^{-1}tk \in T. $$
  Two elements of $T$ are $K$-conjugate if and only 
  if they are $N(T)$-conjugate. So the equivalences continue:
  \begin{eqnarray*}
   &\Iff& \exists w\in N(T),\ k^{-1}tk = w^{-1}tw \\
   &\Iff& \exists w\in N(T),\ wk^{-1}tkw^{-1} = t \\
   &\Iff& \exists w\in N(T),\ kw^{-1} \in C_K(t) \\
   &\Iff& k \in C_K(t) N(T). 
  \end{eqnarray*}
  This chain of equivalences establishes the first claim.

  For the second claim, let us first note that since $T$ is commutative
  $C_K(t) \geq T$, and since $T$ is connected, $C_K(t)^0 \geq T$. 

  Plainly $C_K(t) N(T) \supseteq C_K(t)^0 N(T)$, so our next task is to
  show $C_K(t) \subseteq C_K(t)^0 N(T)$, which will establish
  $C_K(t) N(T) = C_K(t)^0 N(T)$. Let $c\in C_K(t)$. Then
  $c T c^{-1} \leq c C_K(t)^0 c^{-1} = C_K(t)^0$, so $c T c^{-1}$ is
  another maximal torus of the compact connected group $C_K(t)^0$.
  Hence $\exists d\in C_K(t)^0$ such that $d (c T c^{-1}) d^{-1} = T$.
  So $dc \in N(T)$, and $c \in d^{-1} N(T) \subseteq C_K(t)^0 N(T)$,
  completing this task.

  \newcommand\tildew{{\widetilde w}}
  The components of $C_K(t)^0 N(T)/T$ are the orbits of the connected
  group $C_K(t)^0$ through the discrete set $N(T)/T$. 
  Let $\tildew \in N(T)$ lie over $w\in W$. 
  Then $C_K(t) \tildew T/T = C_K(t)/T w \iso C_K(t)/T$, 
  as claimed.
  The $T$-fixed points on the component $C_K(t)^0 \tildew T/T$ are
  $W_t w$. Two components are equal if and only if their $T$-fixed
  points are the same, so the components are indexed by $W_t \backslash W$.

  We turn to the third claim.
  Let $Z$ denote the identity component of the center of $C_K(t)^0$. 
  So $Z$ is a connected subgroup of $K$ commuting with the maximal torus $T$,
  and hence $Z \subseteq T$.
  We now claim that any point $kT$ stabilized by $t$ is also stabilized by $Z$.
  
  By the first claim, we can factor $k$ as $k = cw$, where $c\in C_K(T)^0,
  w\in N(T)$. Then for any $z\in Z$,
  $$ z k T = z c w T = c z w T = c w t' T = c w T = k T $$
  where since $z\in T$, we have $t' = w^{-1} z w$ is also in $T$.

  So for $t$ to occur as a finite stabilizer, $Z$ must be trivial,
  meaning $C_K(t)^0$ must be semisimple.

  For the converse, we know from lemma \ref{lem:genericstab} that the
  generic $T$-stabilizer on $C_K(t)^0/T$ is just $Z(C_K(t)^0)$. 
  This latter group will be finite if and only if $C_K(t)^0$ is semisimple.
\end{proof}

The center $Z(K)$ supplies dull examples of elements of $T$ with
semisimple centralizer (namely, all of $K$). 
If $K=SU(n)$, then there are no other examples.
For a first taste of what can happen in other Lie types, 
consider the diagonal matrix $t = diag(-1,-1,-1,-1,+1)$ in $SO(5)$, 
which has $C_K(t)^0 = SO(4)$. This element fixes $|W_t\backslash W| = 2$ 
copies of $SO(4)/T$ in $SO(5)/T$.

A conjugacy class in $K$ is called {\dfn special} if the centralizer
of some (hence any) element of the class is semisimple. We will
typically use representatives $t\in T$, which we can do since $T$
intersects every conjugacy class. Since semisimplicity is a Lie
algebra phenomenon, it is enough to check that $\lie{c}_{\lie{k}}(t)$
is semisimple.

To analyze these special conjugacy classes, we run down the
arguments from \cite{BdS}, where greater detail can be found.
Recall that for $K$ simple, the {\dfn affine Dynkin diagram} of $K$ is
formed from the simple roots and the {\em lowest} root, which we'll
denote $\omega$.

\begin{proposition}\label{prop:specialclasses}\cite{BdS}
  Assume $K$ is semisimple, so that there are special conjugacy classes,
  and the universal cover $\widetilde K$ is again compact.
  The special conjugacy classes in $K$ are images of those in $\widetilde K$, 
  so it suffices to find those of $\widetilde K$.

  Now assume $K$ simple.
  The special conjugacy classes in $\widetilde K$ correspond $1:1$ to
  the vertices of the affine Dynkin diagram of $\widetilde K$ (or of $K$).
  To find an element of the special class corresponding to a vertex $v$,
  find an element $t\in T$ annihilated by all the roots
  in the affine diagram {\em other than} $v$.
  These $t\in T$ exist and are special. 

  The simple roots and the lowest root satisfy a unique linear dependence
  $\omega + \sum_{\alpha} c_\alpha \alpha = 0$, which we use to
  define the coefficients $\{c_\alpha\}$.
  The adjoint order (meaning, in $K/Z(K)$) 
  of a special element corresponding to a simple
  root $\alpha$ is $c_\alpha$. 
  (For example, a special element is central 
  if and only if the corresponding $c_\alpha$ coefficient is $1$.)
\end{proposition}

These coefficients $\{c_\alpha\}$ can be found 
in e.g. \cite[p98]{HumphreysRGCG}.
They are all $1$ for $K=SU(n)/Z_n$, with the consequence that
the identity is the only finite stabilizer, and the weight varieties
are all manifolds. They are all $1$ or $2$ for the classical groups
$SU(n), SO(n), U(n,{\mathbb H})$.

We now recall the role of the {\em Weyl alcove} in analyzing the 
conjugacy classes of $\widetilde K$.

Each conjugacy class in $K$ meets $T$, and two elements of $T$ are
conjugate in $K$ only if they're already conjugate by the action of $W$.
So the space of conjugacy classes is 
$T/W = (\lie{t}/\Lambda)/W = \lie{t}/(\Lambda \rtimes W)$,
where $\Lambda$ is
the {\dfn coweight lattice} $\ker (\exp: \lie{t} \to T)$.
If $K = \widetilde K$, then this semidirect product is again a
reflection group, the affine Weyl group $\widehat W$; this is the
reason it's convenient to work with $\widetilde K$.
This group is generated by the reflections in the hyperplanes
$\langle \alpha, \cdot\rangle = 0$ for $\alpha$ simple and
$\langle \omega, \cdot\rangle = -1$. (If $K$ is not simple,
then there are several lowest roots, an uninteresting complication
we will ignore.)

The {\dfn Weyl alcove} $A \subset \lie{t}$ is a fundamental region
for $\widehat W$, defined by $\langle \alpha, \cdot\rangle \geq 0$
for $\alpha$ simple and $\langle \omega, \cdot\rangle \geq -1$.
This, and the analysis above, ensure that for $K=\widetilde K$ the image
$\exp(A) \subseteq T \leq K$ intersects each conjugacy class of $K$
in exactly one point. Since the map $A \to \exp(A)$ is a homeomorphism
we may also sometimes refer to $\exp(A)$ as the Weyl alcove.

\begin{corollary}
  If $K$ is simply connected, then the finite stabilizers in the
  action of $T$ on $K/T$ are the Weyl conjugates of the vertices of
  the Weyl alcove $\exp(A)$.
  If $K$ is not simply connected, then the finite stabilizers are the
  images of the finite stabilizers from the universal cover $\widetilde K$.

  If $K$ is a centerless classical group, then $\Gamma$ is contained in
  the $2$-torsion subgroup of $T$. If $K$ is classical but not centerless,
  then $\Gamma$ is contained in the preimage of the $2$-torsion in 
  the torus of $K/Z(K)$.
\end{corollary}

The Weyl alcove also comes up when working examples, in that the torus
$T$ can be pictured as a quotient of the polytope $\bigcup W\cdot A$
made from the union of the Weyl reflections of the Weyl alcove. 
We will call this the {\dfn Tits polytope} after its relation to the
Tits cone in the corresponding affine Kac-Moody algebra,
and to the finite Tits building living on its surface
(neither of which are relevant here, thankfully).
The Tits polytope tesselates the vector space $\lie{t}$,
under translation by the coweight lattice $\Lambda$.
The exponential map $\bigcup (W\cdot A) \to T$ is onto, but only
$1:1$ on the interior of the polytope.

\subsubsection{Example: $K=G_2$}
This group is both centerless and simply connected, hence the
unique Lie group with its Lie algebra. Its Weyl alcove is a 
$30^\circ$-$60^\circ$-$90^\circ$ triangle. 

\begin{center}
%\raisebox{1.5in}
  {
    \begin{minipage}[t]{0.62\linewidth}
      We picture $\Gamma \iso Z_2 \times Z_6$
      inside the Tits hexagon $\bigcup (W\cdot A) \subset \lie{t}$,
      seen on the right.
      (The Weyl alcove itself is the 
      $30^\circ$-$60^\circ$-$90^\circ$ triangle with
      vertices $\{1,\tau,\theta\}$, of orders $1,2,3$.)
      This hexagon tiles the plane $\lie{t}$
      under translation by the coweight lattice $\Lambda$.
      The twelve black dots are labeled by the elements of $\Gamma$
      they exponentiate to, some of which occur again as white dots.
      The element $\sigma$ is a Weyl conjugate of $\tau$.
    \end{minipage} }
  \hfill
  \raisebox{-1.2in}
  {
    \begin{minipage}[t]{0.30\linewidth}
      \epsfig{file=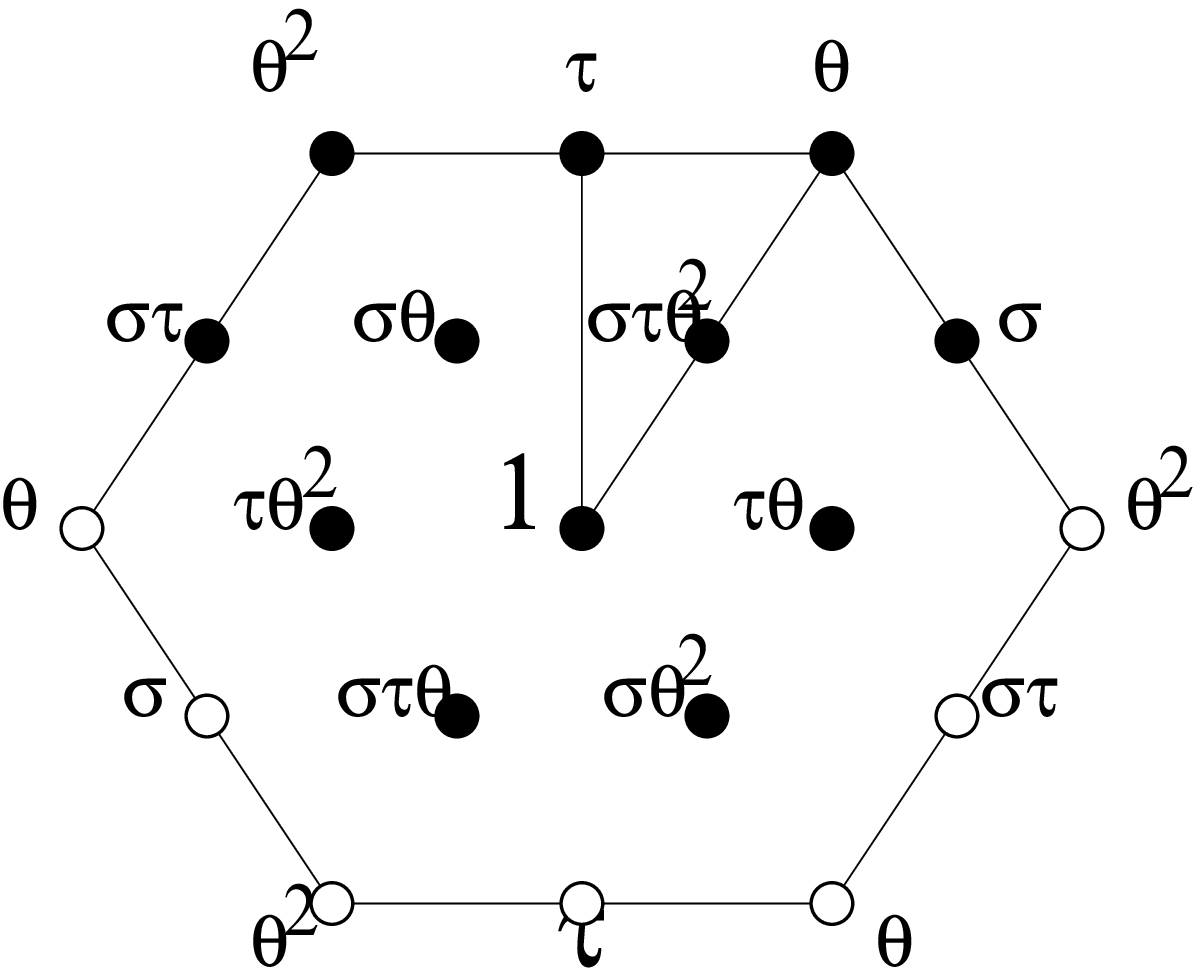,height=1.5in}
    \end{minipage} }
\end{center}

The centralizers of $\tau$ and $\theta$ are $SO(4)$ and $SU(3)$ respectively.
Hence the fixed point sets for the twelve elements of $\Gamma$ are
\begin{itemize}
\item all of $X$, for $w=1$,
\item two $SU(3)/T$s, for $\theta$ and $\theta^2$, 
\item three $SO(4)/T = (\P^1)^2$s, for $\tau,\sigma,\sigma\tau$, and
\item six $U(2)/T = \P^1$s, for each of the other six elements in the interior
  of the hexagon.
\end{itemize}
Only the first six, which lie on {\em corners} of permuted Weyl alcoves, 
are actually finite stabilizers. The last six lie only on edges.

\subsection{The inertial cohomology groups $NH^{*,t}_T(K/T)$}

For $t\in T$, we've already computed the fixed point set as the
disjoint union $\bigcup_{w \in W_t\backslash W} C_K(t)^0 wT/T$.
Therefore
$$ NH^{*,t}_T(K/T) = \bigoplus_{w \in W_t\backslash W} H^*_T(C_K(t)^0 wT/T).$$
The moment polytope of the component $C_K(t)^0 wT/T$ is the convex hull
of the points $W_t w \cdot \lambda$. 

To compute the Weyl group $W_t$ of $C_K(t)^0$, we only need
a set of reflections that generates $W_t$. These are the reflections 
in the roots perpendicular to the (possibly internal) walls of the
Tits polytope passing through the point $t$. For there to be enough
to make $t$ a finite stabilizer, $t$ has to lie on a vertex of a 
permuted Weyl alcove (as already proven in Proposition
\ref{prop:specialclasses}).

The cohomology $H^*_T(C_K(t)^0 wT/T)$ of one component has a basis
given by the equivariant classes of Schubert varieties. For later
purposes we will also be interested in the classes of permuted 
Schubert varieties, for which our reference is \cite{Go:thesis}.

\subsection{The product structure on $NH^{*,\diamond}_T(K/T)$}

Let $t,s \in T$. Then since
$NH^{*,t}_T(K/T)$ is the direct sum 
$\bigoplus_{w \in W_t\backslash W} H^*_T(C_K(t)^0 wT/T)$
and 
$NH^{*,s}_T(K/T)$ is a similar direct sum, to understand their product
it is enough to consider the product from two summands, 
$H^*_T(C_K(t)^0 wT/T) \times H^*_T(C_K(s)^0 vT/T)$.
The definition of the $\smile$ product requires us to restrict classes
from $C_K(t)^0 wT/T$ and $C_K(s)^0 vT/T$ to their intersection,
multiply together and by the virtual fundamental class, and then
push into $C_K(ts)^0 N(T)/T$.

% For $K$ a compact connected Lie group with maximal torus $T$, there
% is a canonical factorization of $T$ as a central product of two
% subtori. The first is the identity component $Z$ of the center of $K$.
% The second $S$ is the identity component of the intersection of
% $T$ with the commutator subgroup of $K$. One can check on the Lie algebra
% level that $T = ZS$, and $Z\cap S$ is discrete.
% Factor $T = ZS$ as explained above, where $K = C_K(t)^0$. 

\begin{lemma}\label{lem:Cst}
  Let $C_K(t,s)$ denote the intersection $C_K(t) \cap C_K(s)$, and
  $C_K(t,s)^0$ its identity component. Let $W_{t,s}$ denote its Weyl group.
  \begin{enumerate}
  \item The intersection $(C_K(t)^0 wT/T)\ \cap\ (C_K(s)^0 vT/T)$
    is fixed pointwise by $t$ and $s$.
    It is a finite union of homogeneous spaces for $C_K(t,s)^0$,  
    namely $\bigcup_{W_{t,s} u} C_K(t,s)^0 uT/T$,
    where the components are indexed by those cosets $W_{t,s} u$
    contained in the intersection $W_t w \cap W_s v$.
  \item The obstruction bundle over $C_K(t,s)^0 uT/T$ is trivial 
    if $K$ is a classical group.
  \end{enumerate}
\end{lemma}

\begin{proof}
  Since $s\in T\leq C_K(t)^0$, we can compute $C_K(s) \cap C_K(t)^0$ as
  the centralizer in $C_K(t)^0$ of $s$. Likewise, the intersection
  $(C_K(t)^0 wT/T)\ \cap\ (C_K(s)^0 N(T)/T)$ can be computed as the 
  $s$-fixed points on $C_K(t)^0 wT/T \iso C_K(t)^0/T$.
  Then apply Lemma \ref{lem:stabs} to the case of $C_K(t)^0$.
  
  Now assume $K$ is classical, and let $g_1,g_2 \in \Gamma$. Then
  $g_1, g_2, (g_1 g_2)^{-1}$ each act on $K/T$ with order $1$ or $2$.
  Hence their three logweights on any line each live in $\{0,1/2\}$,
  and can't add up to $2$. So the obstruction bundle is trivial.
\end{proof}

\begin{corollary}\label{cor:classical}
  Let $K$ be a centerless classical group. Then the product map
  $$ \smile: H^*_T(C_K(t)^0 wT/T) \times H^*_T(C_K(s)^0 vT/T) 
  \to H^*_T(C_K(st)^0 vT/T) $$
  is given by
  $$ \alpha \smile \beta = \sum_{W_{t,s} u \subseteq W_t w \,\cap\, W_s v} 
  (\bar e_3)_*(e_1^*(\alpha) \times e_2^*(\beta)) $$
  where $e_1, e_2, \bar e_3$ are the inclusions of $C_K(s,t)^0 uT/T$
  into $C_K(t)^0 uT/T$, $C_K(s)^0 uT/T$, $C_K(ts)^0 uT/T$, respectively.
\end{corollary}

The maps $e_i^*$ and $(e_i)_*$ between the non-equivariant
cohomologies of these homogeneous spaces have been studied in
\cite{Purbhoo}, in part for the application in \cite{BerensteinSjamaar}
to asymptotic branching rules. 

\subsubsection{Example: $K=SO(5)$} The Weyl alcove of $\widetilde K = Spin(5)$
is a $45^\circ$-$45^\circ$-$90^\circ$ triangle, and the group $\Gamma$
in $\widetilde K$ is exactly the $2$-torsion in $T$, all of whose
elements are finite stabilizers. Its quotient in $SO(5)$ is the
two-element group $\{{\bf 1}, t := diag(-1,-1,-1,-1,+1)\}$.

Hence, there are only two summands in $NH_T^{*,\Gamma}(SO(5)/T)$.
By Remark \ref{rem:idprod}, 
the only difficult product is from the $t$ summand, squared, back to
the identity summand. In this case $e_1$ and $e_2$ are the identity,
so the only map of interest is 
$(\bar e_3)_* : H^*_T((SO(5)/T)^t) \to H^*_T(SO(5)/T)$. This is perhaps
best computed via the techniques in the last section.

\subsection{The kernel of the inertial Kirwan map}

Finally, we need to compute the kernel of the map from
$NH^{*,\Gamma}_T(K\cdot\lambda)$ to $H^*_{CR}(K\cdot\lambda //_\mu T)$.
Breaking this up by $t\in\Gamma$, 
and then into components of $(K\cdot\lambda)^t$, this kernel is 
the direct sum of the kernels of each of the ordinary Kirwan maps
$$ H^*_T(C_K(t)^0 u \cdot \lambda) \to 
H^*(C_K(t)^0 u \cdot \lambda //_\mu T).$$
This kernel is computed in \cite{Go:thesis}; it is spanned by the
classes of those permuted Schubert varieties whose image under the
moment map misses $\mu$.

%% 9. EGs: Toric varieties -- Tara
%%     EGtoric.tex
%%        * \Z coefficients
%%        * Relation with Borisov-Chen-Smith, etc
%%        * Relation with Lerman-Tolman

\section{Toric varieties}
\label{sec:toricvarieties}
In this section, we use our results to compute the Chen-Ruan
cohomology of certain toric orbifolds.  We first discuss the
symplecto-geometric construction of toric orbifolds, as described
by Lerman and Tolman \cite{LT:toricorbifolds}.  We remark on the
coefficients in the toric case and we compute an example in full
detail. Finally, we relate our results to those of Borisov, Chen
and Smith \cite{BCS:toricvarieties}.

\subsection{Symplectic toric orbifolds}\label{subse:toric}

In \cite{LT:toricorbifolds}, Lerman and Tolman study Hamiltonian
torus actions on symplectic orbifolds, and define and classify
symplectic toric orbifolds.   They consider the case where the orbifold
is {\bf reduced}: there is no global stabilizer.  These reduced toric orbifolds,
then, are in one-to-one
correspondence with labeled simple rational polytopes.

\begin{definition}
Let $\algt$ be a $d$-dimensional vector space with a distinguished
lattice $\ell$; let $\algt^*$ be the dual space.  A convex
polytope $\Delta\subset \algt^*$ is {\dfn rational} if
$$
\Delta = \bigcap_{i=1}^N \{ \alpha\in\algt^* \ |\  \langle \alpha,
y_i\rangle \geq \eta_i \}
$$
for some $y_i\in \ell$ and $\eta_i\in\Q$.  A {\dfn facet} is a
face of codimension $1$.  A $d$-dimensional polytope is {\dfn
simple} if exactly $d$ facets meet at every vertex.  A {\dfn
labeled polytope} is a convex rational simple polytope along with
a positive integer labeling each facet.
\end{definition}

To establish a one-to-one correspondence between symplectic toric
orbifolds and labeled polytopes, Lerman and Tolman mimic Delzant's
construction of (smooth) toric varieties as symplectic quotients.
The labeled polytope $\Delta$ is uniquely described as the
intersection of half-spaces,
$$
\Delta = \bigcap_{i=1}^N\ \{ \alpha\in\algt^* \ |\  \langle \alpha,
y_i\rangle \geq \eta_i \},
$$
where $N$ is the number of facets, the vector $y_i\in \ell$ is the
primitive inward-pointing normal vector to the $i$th facet, and
$m_i$ is the positive integer labeling that facet. Define the map
$\varpi : \R^N \rightarrow \algt$ defined by sending the $i$th
standard basis vector $e_i$ to $m_i y_i$.  This yields a short
exact sequence,
\begin{equation}\label{eq:tseq}
\xymatrix{
 0 \ar[r] &  \algk \ar[r]^{j} & \R^N \ar[r]^{\varpi} &  \algt
 \ar[r] & 0,
 }
\end{equation}
and  its dual
\begin{equation} \xymatrix{
 0 \ar[r] &  \algt^* \ar[r]^{\varpi^*} & (\R^N)^* \ar[r]^{j^*} &
 \algk^* \ar[r] & 0,
 }
\end{equation}
where $\algk = {\mathrm ker}(\varpi).$  Let $\T^N = \R^N / \Z^N$
and $T = \algt/\ell$, and let $K$ denote the kernel of the map
$\T^N \to T$ induced by $\varpi$.  Then the Lie algebra of $K$ is
$\algk$. Tolman and Lerman then prove that $Y_\Delta =  \C^N/\! /
K$ is the unique symplectic toric orbifold with moment polytope
$\Delta$.

\begin{figure}[h]
\begin{center}
\psfrag{2}{\tiny{2}}
\includegraphics[height=0.1in]{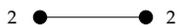}
\end{center}
\caption{The labeled polytope corresponding to $\C^2 /\!
/(S^1\times \Z_2)$. }
 \label{fig:DisconnTorus}
\end{figure}

As noted above, this construction does not include the possibility of a global finite
stabilizer. It may, however, result in a quotient by a disconnected
subgroup $K$ of $\T$. For example, for the polytope shown in
Figure~\ref{fig:DisconnTorus}, using the above construction, we find
that $K\iso S^1\times \Z_2$.  In order to obtain a connected group
$K$, it is sufficient, though not necessary, to assume that the labels
are all $1$.  When $K$ is connected, we may use the techniques
developed in Sections~\ref{sec:starproduct} and \ref{sec:surjectivity}
to compute the Chen-Ruan cohomology of $Y_\Delta$.  When $K$ is 
disconnected, we need an additional argument to see that the surjectivity
remains true (over $\Q$).

\begin{proposition}
  Let $Y_\Delta= \C^N/\! /K$ be a symplectic toric orbifold.  Then
  there is a surjection
  $$   H_K^*(\C^N;\Q) \longrightarrow H^*(Y_\Delta;\Q).   $$
\end{proposition}

\begin{proof}
  We can write $K= T\times G$, where $T$ is a connected torus, and $G$ is a
  finite abelian group (possibly trivial). Since $G$ is finite, we may identify
  $$
  H_{T\times G}^*(\C^N;\Q)\iso H_T^*(\C^N/G;\Q).
  $$
  Similarly, for a level set $\Phi^{-1}(0)$, there is an isomorphism
  $H_{T\times G}^*(\Phi^{-1}(0);\Q)\iso H_T^*(\Phi^{-1}(0)/G;\Q)$.
  Now, $\C^N/G$ is a symplectic orbifold, and $T$ acts on $\C^N/G$ in
  a Hamiltonian fashion.  Moreover, the level set for a $T$ moment map
  may be identified as $\Phi^{-1}(0)/G$, the $K$-level set modulo $G$.
  One may extend Kirwan surjectivity to orbifolds using the techniques
  detailed in \cite{LMTW} and thus conclude that
  $$
  H_T^*(\C^N/G;\Q) \longrightarrow H^*_T(\Phi^{-1}(0)/G;\Q).
  $$
  Finally, we identify $H^*_T(\Phi^{-1}(0)/G;\Q)\iso H^*(Y_\Delta;\Q)$, 
  completing the proof.
\end{proof}

From the labeled moment polytope, it is not difficult to determine the
stabilizer group of any point.  There are two contributions: a
contribution that can be seen in the underlying polytope, and a
contribution from the facet labels. Given a point in $Y_\Delta$, let
$f$ be the minimal face in $\Delta$ that contains its image under the
moment map, let $V$ be the subspace containing $f$, and $V_{\Z} =
V\cap \Z^N$.  Let $E$ be the set of primitive vectors incident to a
vertex of $f$ that are not in the subspace $V$. Let $\pi$ be the
projection from $\Z^N$ to $\Z^N/V_{\Z}$. The contribution to the
stabilizer coming from the polytope is the group
$$
\Gamma^1_f = (\Z^N/V_{\Z})\big/ \Span_{\Z}(\pi(E)).
$$
Now we let $F_i$ denote the facet of $\Delta$ with label $\ell_i$.
The contribution to the stabilizer from the labels is
$$
\Gamma^2_f = \gcd ( \ell_i\ |\ f\subseteq F_i ).
$$
And so the stabilizer to the point we are considering is
\begin{equation}\label{eq:stab}
\Gamma_f = \Gamma^1_f \times \Gamma^2_f.
\end{equation}

\subsection{A comment on coefficients}\label{se:integers}

The surjectivity statement in Theorem~\ref{thm:surjectivity}
follows directly from Kirwan's Theorem~\ref{thm:kirwan}.  As a
consequence, Theorem~\ref{thm:surjectivity} holds over $\Z$
whenever Theorem~\ref{thm:kirwan} holds over $\Z$ for each of the orbistrata.
Even for a smooth toric variety, however, it is not immediately clear
(though it is true \cite{HT}) that
Kirwan's proof of surjectivity generalizes to integer
coefficients. Here we first establish Kirwan's result over $\Z$ for
weighted projective spaces. 

\begin{proposition}\label{prop:Zsurj}
Suppose that $S^1$ acts on $\C^N$ linearly with positive weights.
Then $\C^N /\! /S^1$ is a {\bf weighted projective space}, and the ring
homomorphism
\begin{equation}
     \kappa\ :\ H_{S^1}^*(\C^{N};\Z) \to H^*(\C^N/\! /S^1; \Z)
\end{equation}
is a surjection.
\end{proposition}

\begin{proof}
Since $S^1$ acts linearly on $\C^N$ with positive weights, the moment
map is of the form
$$
\Phi(z_1,\dots,z_N) = b_1 |z_1|^2 + \cdots + b_N |z_N|^2 + C,
$$
where $b_i\geq 0$ are positive multiples of the positive weights of
the circle action, and $C$ is a constant.  Hence, a level set at a
regular value is homeomorphic to a sphere and the resulting quotient
$S^{2N-1}/S^1$ is precisely a weighted projective space.  Indeed,
every weighted projective space arises in this way.

Surjectivity holds over $\Z$ in these examples by results similar to
Propositions~7.3 and 7.4 in \cite{TW:symplecticquotients}, where
Tolman and Weitsman show that surjectivity over $\Z$ depends only on
the integral cohomology of the circle fixed points being
torsion-free. Although Tolman and Weitsman's results require the
original manifold to be compact, because we have
\begin{enumerate}
\item $\Phi$ is proper and bounded below; and
\item $\C^N$ has only finitely many fixed point components,
\end{enumerate}
we may generalize their propositions to this setting.   Thus, for
surjectivity to hold over $\Z$, it is sufficient for the integral
cohomology of each fixed point component to be torsion-free. Now, in
this case, the only circle  fixed point is the origin in $\C^N$, and
of course a  point has torsion-free integral cohomology.  This
completes the proof.
\end{proof}

In a preliminary version of this paper, we stated a stronger (but
incorrect) version of Proposition~\ref{prop:Zsurj}.  It applied to a
class of toric orbifolds with labeled polytope of a certain
combinatorial type, including, but not limited to, weighted projective
spaces. To prove it, we established that the critical sets of
$||\Phi||^2$ have torsion-free integral cohomology.  This condition is
in fact not sufficient to deduce surjectivity over $\Z$.  Indeed, the
manifold $X$ in \cite[Lemma~7.1]{TW:symplecticquotients} must be {\em
pointwise} fixed by the torus $T$, so we may only apply this Lemma to
the $T$-fixed points of the original manifold; however, unless $T$ is
a circle, there are critical sets $C$ of $||\Phi||^2$ that are only
fixed by a subtorus $K$ of $T$.  As noted in
Remark~\ref{rem:surjinteger}, to prove surjectivity over $\Z$, we need
not only that $H^*(C;\Z)$ is torsion free, but also that
$H_K^*(C;\Z)$ is torsion free.

Theorem~\ref{thm:kirwan} does hold over $\Z$ more generally, though
not for all toric orbifolds.  For example, consider the product of
weighted projective spaces $X = \P_{1,1,2} \times \P_{1,1,2}$.  Using
the K\"unneth theorem over $\Z$, it is straight forward to compute the
cohomology of this space.  Because $\P_{1,1,2}$ has $2$-torsion in its
cohomology, the $\mathrm{Tor}$ terms do not vanish, and as a result,
the cohomology of the product has $2$-torsion in certain odd
degrees. As a result, surjectivity cannot possibly hold over $\Z$
because the equivariant cohomology of affine space has cohomology only
in even degrees. The failure in this example is entirely due to the
fact that there is a point with stabilizer $\Z_2\times\Z_2$.  For a
toric orbifold whose polytope has labels all equal to $1$, to ensure
surjectivity over $\Z$, it is at the very least necessary to assume
that the isotropy at each point of the orbifold is a cyclic
group. This is being more closely investigated by Tolman and the
second author \cite{HT}.

We now use Proposition~\ref{prop:Zsurj} to deduce an integral
version of Theorem~\ref{thm:surjectivity} for a weighted
projective space $\C^N/\! /S^1$.  Each orbistratum is
itself a circle reduction of some $\C^k\subseteq \C^N$ and thus also a
weighted projective space, immediately
implying surjectivity over $\Z$ for each piece of the inertial cohomology.
Hence we have the following corollary.

\begin{corollary}
Suppose that $S^1$ acts on $\C^N$ linearly with positive weights.
Then the ring homomorphism
\begin{equation}
     \kappa_{NH}\ :\ \srous(\C^{N};\Z) \to H^*_{CR}(\C^N/\! /S^1; \Z)
\end{equation}
is a surjection.
\end{corollary}

Finally, we note that whenever the surjectivity result holds over
$\Z$ for each orbistratum, the kernel computations of
Theorem~\ref{th:noncompactTW} extend to inertial cohomology with $\Z$
coefficients.  As the examples presented below are all weighted
projective spaces or smooth toric varieties, surjectivity and the
kernel computations do hold over $\Z$.

\subsection{The Chen-Ruan cohomology of a weighted projective space}

We now present an example to demonstrate the ease of computation
of inertial cohomology and of the kernel of the surjection to
Chen-Ruan cohomology of the reduction.

\begin{example}\label{eg:P123}\rm
Let $\Delta$ be the moment polytope in
Figure~\ref{fig:weightedP2}.

\begin{figure}[h]
\begin{center}
\psfrag{1}{\tiny{1}}
\includegraphics[height=1in]{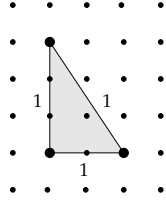}
\end{center}
\caption{The labeled polytope corresponding to $\big(\C_{(1)}\oplus
\C_{(2)}\oplus \C_{(3)}\big)/\! / S^1$.}
 \label{fig:weightedP2}
\end{figure}

\noindent In this case, if we let $Y = \C^3$, and follow the above
construction, then $K\cong S^1$ acts on $Y$ by
$$
t\cdot (z_1, z_2, z_3) = (t\cdot z_1, t^2\cdot z_2, t^3\cdot z_3).
$$
Then as an $S^1$ representation, $Y = \C_{(1)}\oplus
\C_{(2)}\oplus \C_{(3)}$. This action is Hamiltonian, with moment
map
$$
\Phi (z_1, z_2, z_3) = |z_1|^2 + 2 |z_2|^2 + 3|z_3|^2.
$$
Any positive real number is a regular value of $\Phi$, and the
symplectic reduction is a weighted projective space $Y_\Delta =
Y/\! / S^1 = {\P}^2_{1,2,3}$.  Changing the regular value at which we
reduce only changes the symplectic form on $Y_\Delta$. We now
compute the Chen-Ruan cohomology of the symplectic reduction, by
computing $NH_{S^1}^{*,\Gamma}(Y)$ and computing the kernel of the
surjection $\kappa : NH_{S^1}^{*,\Gamma}(Y)\rightarrow
H_{CR}^*(Y_\Delta)$.  In this case, the surjection actually holds
over $\Z$, so we will assume integer coefficients for the
remainder of the example.

First, we notice that the finite stabilizers for the $S^1$ action
are the square and cube roots of unity inside $S^1$.  Thus, the
group that they generate is the set
$\Gamma = \{ \zeta_k = \exp(2\pi i k/6)\ |\ k =
0,\dots,5\}\iso \Z_6$ of sixth roots of unity.
For each $\zeta_i$, $Y^{\zeta_i}$ is
contractible, so $NH_{S^1}^{*,\zeta_i}(Y)$ is free of rank $1$.
To compute $NH_{S^1}^{*,\Gamma}(Y)$, we now refer to the
following table.
\renewcommand{\arraystretch}{1.3}
\begin{equation}\label{eq:123table}
\begin{array}{c|c|c|c|c|c|c|}
g & \zeta_0 & \zeta_1 & \zeta_2 & \zeta_3 & \zeta_4 & \zeta_5 \\
\hline Y^g & Y & \{ 0\} & \C_{(3)} & \C_{(2)} & \C_{(3)} &
\{ 0\} \\ \hline \Year_{\C_{(1)}}(g) & 0 & \frac{1}{6} & \frac{1}{3} & \frac{1}{2} & \frac{2}{3} & \frac{5}{6}  \\
\hline \Year_{\C_{(2)}}(g) & 0 & \frac{1}{3} & \frac{2}{3}  & 0 &
\frac{1}{3} & \frac{2}{3} \\ \hline \Year_{\C_{(3)}}(g)  & 0 &
\frac{1}{2}  & 0 & \frac{1}{2}  & 0 & \frac{1}{2} \\ \hline 2\
\age(g) & 0 & 2& 2 & 2& 2& 4\\ \hline
\genfrac{}{}{0pt}{0}{\mbox{generator of}}{NH_{S^1}^{*,g}(Y)} & 1 &
\alpha & \beta & \gamma & \delta & \eta\\ \hline
\end{array}
\end{equation}
Thus, $NH_{S^1}^{*,\Gamma}(Y)$ is generated as a free module
over $H_{S^1}^*(pt)$  by the elements $1$,
$\alpha$, $\beta$, $\gamma$, $\delta$, and $\eta$.
Let $u$ denote the degree 2 generator of
$H_{S^1}^*(pt)$.  We determine the product
structure, by computing all pairwise products of these
generators.

\begin{remark}
By Corollary~\ref{co:productsthesame} we may abuse notation and write
$\star$ for the product on $NH_{S^1}^{\ast,\Gamma}(Y)$. 
\end{remark}

We demonstrate the calculation of $\eta\star\eta$
using \eqref{eq:starproduct}, and then give the multiplication
table.  We know that $\eta\star\eta\in NH^{*,\zeta_4}_{S^1}(Y)$,
so we must compute what multiple of $\delta$ it is. Since $\eta$
is supported on the $(g=\zeta_5)$-piece of $Y$, we must compute
the logweights of $\zeta_5$ and $\zeta_4$ on $\nu Y^g$, the normal
bundle to the fixed point set $Y^g=\{0\}$. Here $\nu Y^g =
\C_{(1)}\oplus\C_{(2)}\oplus \C_{(3)}$. We find that
\begin{eqnarray*}
\Year^{\{ 0\} }_{(1)}(\xi_5)+\Year^{ \{
0\} }_{(1)}(\xi_5) -\Year^{ \{ 0\} }_{(1)}(\xi_4) & = & \frac{5}{6} + \frac{5}{6}-\frac{2}{3}=1,\\
\Year^{\{ 0\} }_{(2)}(\xi_5)+\Year^{ \{
0\} }_{(2)}(\xi_5) -\Year^{ \{ 0\} }_{(2)}(\xi_4) & = & \frac{2}{3} + \frac{2}{3} - \frac{1}{3} = 1, \mbox{ and}\\
\Year^{\{ 0\} }_{(3)}(\xi_5)+\Year^{ \{ 0\} }_{(3)}(\xi_5)
-\Year^{ \{ 0\} }_{(3)}(\xi_4) & = & \frac{1}{2} + \frac{1}{2} - 0
= 1.
\end{eqnarray*}
Using \eqref{eq:starproduct}, we need only calculate
$(\eta\star\eta)|_{\{0\}}$, where the product is
$$
\eta|_{\{0\}}\cdot\eta|_{\{0\}}
\prod_{i=1}^3e(\C_{(i)})^{\Year^{\{ 0\} }_{(i)}(\xi_5)+\Year^{ \{
0\} }_{(i)}(\xi_5) -\Year^{ \{ 0\} }_{(i)}(\xi_4)} =1\cdot 1\cdot
(u)^1\cdot (2u)^1\cdot (3u)^1
$$
for $\{0\}$ the fixed point of $Y^{\zeta_5\zeta_5}=Y^{\zeta_4}$.
Thus $\eta\star\eta= 6u^3\delta$. Similarly we compute the other
products:
\begin{equation}\label{eq:multtable}
\begin{array}{c||c|c|c|c|c|}
\star & \alpha & \beta & \gamma & \delta & \eta\\ \hline \hline
\alpha & 3u\beta & 2u\gamma & 3u\delta & \eta & 6u^3 \\ \hline
\beta & & 2u\delta & \eta & 2u^2 & 2u^2\alpha \\ \hline \gamma & &
& 3u^2 & u\alpha & 3u^2\beta
\\ \hline \delta & & & & u\beta & 2u^2\gamma \\ \hline \eta & & & & &  6u^3\delta \\ \hline
\end{array}.
\end{equation}
Thus, as a ring,
$$
NH_{S^1}^{*,\Gamma}(Y;\Z) \iso \Z[u,\alpha , \beta , \gamma , \delta ,\eta]/I,
$$
where $I$ is the ideal generated by the product relations \eqref{eq:multtable}.

Finally, to compute $H_{CR}^*(Y/\! / S^1;\Z)$, we compute the
kernel of the Kirwan map.  Following \cite{K:quotients}, this
kernel is
\begin{equation}\label{eq:kernel}
\ker(\kappa) = \langle 6u^3,\alpha,3u\beta,2u\gamma,3u\delta,\eta \rangle.
\end{equation}
Thus,
$$
H_{CR}^*(Y_\Delta;\Z) \iso \Z[u,\alpha , \beta , \gamma , \delta
,\eta]/J,
$$
where $J$ is the ideal generated by the relations from
\eqref{eq:multtable} and from \eqref{eq:kernel}.

It is sometimes more convenient to compute with the product
$\smile$.  We conclude this example with the computation of $\eta
\smile \eta$. Using \eqref{eq:123table}, we note that
$e_1^*(\eta)\cdot e_2^*(\eta)$ is in $H_T^*(Y^{\zeta_5,\zeta_5})$.
We notice that $Y^{\zeta_5,\zeta_5} = \{ 0\}$, and in the
cohomology of $Y^{\zeta_5,\zeta_5}$, $e_1^*(\eta)\cdot
e_2^*(\eta)=1$. Thus, $\nu Y^{\zeta_5,\zeta_5} = \C_{(1)}\oplus
\C_{(2)} \oplus\C_{(3)}$.  We now check to see which of these
lines are in the obstruction bundle by computing
\begin{eqnarray*}
\Year^{ \{ 0\} }_{ (1)}(\xi_5)+\Year^{ \{ 0\} }_{
(1)}(\xi_5)+\Year^{ \{ 0\} }_{ (1)}(\xi_4) & = & \frac{5}{6} +
\frac{5}{6} + \frac{1}{3} = 2,\\
\Year^{ \{ 0\} }_{ (2)}(\xi_5)+\Year^{ \{ 0\} }_{
(2)}(\xi_5)+\Year^{ \{ 0\} }_{ (2)}(\xi_4) & = & \frac{2}{3} + \frac{2}{3} + \frac{2}{3} = 2, \mbox{ and}\\
\Year^{ \{ 0\} }_{ (3)}(\xi_5)+\Year^{ \{ 0\} }_{
(3)}(\xi_5)+\Year^{ \{ 0\} }_{ (3)}(\xi_4) & = & \frac{1}{2} +
\frac{1}{2} + 0 = 1.
\end{eqnarray*}
Thus, the obstruction bundle is $E|_{Y^{\zeta_5,\zeta_5}} =
\C_{(1)}\oplus\C_{(2)}$, and so the virtual class in this case is
$\varepsilon = 2u^2$.  Finally, we note that the pushforward
$(\overline{e}_3)_*$ will multiply the class $e_1^*(\eta)\cdot
e_2^*(\eta)\cdot\varepsilon$ by $3u\delta$, which is the Euler
class of the normal bundle to $Y^{\zeta_2}$.  Thus,
$\eta\smile\eta = 6u^3\delta$. The other pairwise products can be
easily computed in the same fashion, yielding, of course, the same
multiplication table as in \eqref{eq:multtable}.
\end{example}

\subsection{Relation to toric Deligne-Mumford stacks}

Borisov, Chen and Smith \cite{BCS:toricvarieties} compute the rational
Chen-Ruan Chow ring of a toric Deligne-Mumford stack. This stack is
determined by a combinatorial object called a stacky fan
$\mathbf{\Sigma}$. To each labeled polytope $\Delta$, we may associate
such a stacky fan.  If the labels of the polytope are all $1$, then
the associated fan is the {\bf canonical stacky fan} of
\cite{BCS:toricvarieties}. Just as for smooth toric varieties,
however, there are stacky fans that cannot be defined via labeled
polytopes.  Here we show that the \cite{BCS:toricvarieties} description
of the Chow ring of these stacky fans corresponds to our description,
over $\Q$. While our results are less general in that we cannot
describe the ring for all stacky fans, we have the advantage
that in some cases we can do computations over $\Z$. When $Y$ is a
symplectic toric orbifold, our results agree.

\begin{theorem}\label{thm:relnToBCS}
Let $Y = \C^N/\! / K$ be a toric orbifold that can be realized as
a Deligne-Mumford stack $\mathcal{X}(\mathbf{\Sigma})$. Then there
is a ring isomorphism
\begin{equation}
H_{CR}^*(Y)\iso\frac{\krous(\C^N;\Q)}{\ker\kappa_{NH}} \longrightarrow
A_{CR}^{*} \bigl( \mathcal{X}(\mathbf{\Sigma}) \bigr)
\end{equation}
that divides all degrees in half.
\end{theorem}

\begin{proof}
We first construct the isomorphism of the modules.  In what
follows, all coefficients are taken to be $\Q$.  As in
Section~\ref{subse:toric}, consider the short exact sequence of
groups
$$
1 \to K \to \T^N \to T \to 1
$$
and the corresponding short exact sequence of dual Lie algebras
$$
\xymatrix{
 0 \ar[r] &  \algt^* \ar[r]^{\varpi^*} & (\R^N)^* \ar[r]^{j^*} &
 \algk^* \ar[r] & 0.
 }
$$
Choose $\Phi : \C^N \to \algk^*$ a moment map with $Y =
\Phi^{-1}(0) / K$.

The inertia stack defined in \cite{BCS:toricvarieties} is given as a
quotient of a space $\locfree$ given as the complement of the subspace
arrangement.
Our moment map level set is homotopy equivalent to their $\locfree$.
Indeed, $\locfree$ is endowed with an algebraic $K_{\C}\iso (\C^*)^k$
action, and the $\R^k$ acts on our moment map level set by scaling
inside $\locfree$.

Borisov, Chen and Smith then show that the inertia stack is
$$
\coprod_{g\in \gamma} [ \locfree^g/K_{\C} ],
$$
where the disjoint union is taken over the set $\gamma$ of finite
stabilizers $g\in K_{\C}$.  As such they give a module isomorphism
$$
A_{CR}^{*} \bigl( \mathcal{X}(\mathbf{\Sigma}) \bigr) \iso A^*
\left( \coprod_{g\in \gamma} [ \locfree^g/K_{\C} ] \right).
$$
Each piece $[\locfree^g/K_{\C}]$ of the inertia stack is again a toric
orbifold. We now construct an isomorphism
$$
NH_K^{*,g}(\C^N)\iso\frac{H_K^*((\C^N)^g)}{\ker\kappa} \longrightarrow
A^*([\locfree^g/K_{\C}])
$$
that divides degrees in half.

Borisov, Chen and Smith describe $A^*([\locfree^g/K_{\C}])$ as in
Danilov \cite{danilov}, also referred to as the Stanley-Reisner description.  
For notational simplicity, we explain the
case when $g= id$.  The remaining pieces of the module isomorphism
are derived in the same way, using subgroups of $T$, $\T$ and $K$
and subpolytopes of $\Delta$ as appropriate. We know
$[\locfree/K_{\C}]$ is a $T$-toric variety with moment polytope
$\Delta$. Danilov proved that
$$
A^*([\locfree/K_{\C}]) \iso \frac{\Q[x_1,\dots,x_N]}{(\I,\J)},
$$
where $\deg(x_i) = 1$; $\I$ is the ideal generated by $\prod_{i\in
I}x_i$ for all $I\subset \{ 1,\dots,N\}$ such the the $I$ facets
do not intersect in $\Delta$; and $\J = \{ \sum \alpha_i x_i\ |\
\alpha\in\varpi^*(\algt)^*\}$.

On the other hand, we describe ${H_K^*(\C^N)}/{\ker\kappa}$
following \cite{TW:symplecticquotients}. There is a commutative
diagram
$$
\begin{array}{c}
\xymatrix{ H_{\T}^*\left(\C^N\right) \ar[r]^{r_K^*}
\ar@{->>}[d]_{\kappa^{\T}} & H_K^*\left(\C^N\right)
\ar@{->>}[d]^{\kappa} \\
H_{\T}^*(\Phi^{-1}(0)) \ar[r]_{r_K^*} & H_K^*(\Phi^{-1}(0))
 }
  \end{array}.
$$
The map $\kappa^{\T}$ is a surjection due to an equivariant
version of Theorem~\ref{thm:kirwan}.  The ring $H_{\T}^*(\C^N)$ is
isomorphic to $\Q[x_1,\dots,x_N]$, where $\deg(x_i)=2$.  The
kernel of the restriction map $r_K^*$ is the ideal $\J$.  As
Tolman and Weitsman show explicitly in the proof of
\cite[Theorem~7]{TW:symplecticquotients}, the kernel of $\kappa$
is $\I$. Finally, since the $K$ action is locally free, we have an
isomorphism $H_K^*(\Phi^{-1}(0)) \iso H^*([\Phi^{-1}(0)/K])$. Thus,
identifying the $x_i$'s gives an isomorphism
$$
H^*([\Phi^{-1}(0)/K]) \longrightarrow A^*([\locfree/K_{\C}])
$$
that divides degrees in half.  Putting these together produces a
module isomorphism
$$
\frac{\krous(\C^N;\Q)}{\ker\kappa_{NH}} \longrightarrow
A_{CR}^{*} \bigl( \mathcal{X}(\mathbf{\Sigma}) \bigr).
$$

More generally, we have a commutative diagram in inertial
cohomology
$$
\begin{array}{c}
\xymatrix{ \Trous(\C^N) \ar[r] \ar@{->>}[d]_{\kappa^{\T}_{NH}} &
\krous(\C^N)
\ar@{->>}[d]^{\kappa_{NH}}  &  \kfinitepreorb(\C^N) \ar@{_(->}[l] \ar@{->>}[dl]\\
\Trous(\Phi^{-1}(0)) \ar[r] & \krous(\Phi^{-1}(0)) & \\
NH_{\mathbb{T}}^{\ast,\gamma}(\Phi^{-1}(0)) \ar@{^(->}[u] \ar@{->>}[ur]& &
 }
 \end{array}.
$$
\cite{BCS:toricvarieties} start at $H_{\T}^*(\C^N)$, first restrict
to the level set and then restrict to the subgroup $K$.  As a
result, $NH_{\mathbb{T}}^{\ast,\gamma}(\Phi^{-1}(0))\iso \Q[x_1,\dots,x_N]/\I$ 
is precisely the
numerator $\mathbb{Q}[N]^{\mathbf{\Sigma}}$ in their
\cite[Theorem~1.1]{BCS:toricvarieties}, once degrees are divided
in half. 

%\todo{I think that I simplified the previous section.  Before there
%was some confusion over $\locfree$ and the level set $\Phi^{-1}(0)$.
%Is it more clear now??  Or more confusing now??  As for mentioning the
%obstruction bundles, I don't really think it's necessary.} 
That the module isomorphism is a ring isomorphism is a simple
exercise in comparing the ring structures constructed here and in
\cite{BCS:toricvarieties}.  We leave the details to the reader.
\end{proof}

Finally, we mention the relationship with crepant resolutions. There
is a general conjecture that the Chen-Ruan cohomology 
of an orbifold 
coincides with the ordinary cohomology (or Chow ring) of a crepant
resolution, if one exists \cite{CR:orbH}. In 
\cite[\S 7]{BCS:toricvarieties}, the authors find
necessary and sufficient conditions for a toric orbifold $Y$ to
have a crepant resolution $\widetilde{Y}$. They describe a flat
family $\mathcal{T}\to\P^1$ of rings such that $\mathcal{T}_0$ is
the Chen-Ruan Chow ring of $Y$, and $\mathcal{T}_\infty$ is the
Chow ring of $\tilde{Y}$. Finally, they computed the Chen-Ruan Chow ring
of the weighted projective space $Y = \P_{1,2,1}$, and the Chow ring
of a crepant resolution $\tilde{Y}$, and show that they are not
isomorphic as rings over $\Q$ (though the conjecture still remains
over $\C$).

\begin{figure}[h]
\begin{center}
\psfrag{1}{\tiny{1}}\psfrag{a}{(a)}\psfrag{b}{(b)}
\includegraphics[height=1.5in]{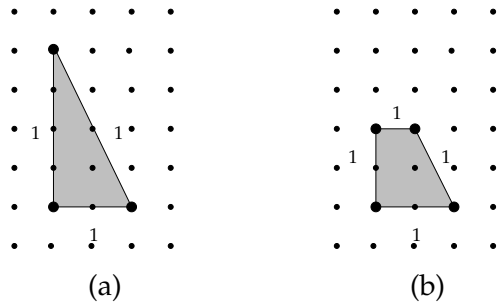}
\end{center}
\caption{
(a) shows the labeled polytope for $Y = \P_{1,2,1}$.  (b) shows the
labeled polytope for $\widetilde{Y}$ its crepant resolution.
}
 \label{fig:crepant}
\end{figure}

We show that even the module structures are not the same over $\Z$.
The space $Y = \P_{1,2,1}$ and its crepant resolution $\widetilde{Y}$
have moment polytopes show in Figure~\ref{fig:crepant}.

\noindent%
Using the methods we have described, we may compute
\begin{equation*}
H_{CR}^*(\P^2_{1,2,1};\Z)  =  \Z[u,\alpha]/\langle u^2-\alpha^2,
2u^3, 2u\alpha\rangle,
\end{equation*}
which yields
\begin{equation*}
H_{CR}^i(\P^2_{1,2,1};\Z)  =  \left\{\begin{array}{ll} \Z & i =
0\\ \Z\oplus\Z & i = 2 \\ \Z\oplus\Z_2 & i = 4 \\ \Z_2\oplus\Z_2 &
i = 2n> 4\\ 0 & \mbox{else}\end{array}\right. .
\end{equation*}
Thus, we see that there is torsion in all higher degrees.  This
certainly does not happen for the crepant resolution, as it is a
smooth variety.  Thus even as modules, the two cohomologies do not
agree over $\Z$.

\end{document}